\documentclass[a4paper]{scrartcl}

\usepackage[english,francais]{babel} 
\usepackage{ucs}
\usepackage[utf8]{inputenc}

\usepackage{graphicx,psfrag,amsfonts,bbm}
\usepackage{amsmath,indentfirst,qsymbols,floatflt}

\RequirePackage{color} 

\RequirePackage{amsmath,amssymb,amsfonts,bbm,latexsym,mathrsfs}

\newcommand{\lf}{\lfloor}
\newcommand{\rf}{\rfloor}

\newcommand{\bc}{{\bf c}}

\newcommand{\ff}{\mathbbm{f}}
\newcommand{\bbf}{{\bf f}}

\newcommand{\bt}{{\bf t}}
\newcommand{\bg}{{\bf g}}

\newcommand{\bx}{{\bf x}}
\newcommand{\bv}{{\bf v}}
\newcommand{\by}{{\bf y}}

\newcommand{\bF}{{\bf F}}

\newcommand{\bT}{{\bf T}}
\renewcommand{\P}{\mathbb{P}}

\newcommand{\bP}{{\bf P}}

\newcommand{\bE}{{\bf E}}

\newcommand{\R}{\mathbb{R}}
\newcommand{\N}{\mathbb{N}}
\newcommand{\F}{\mathbb{F}}
\newcommand{\Z}{\mathbb{Z}}

\newcommand{\T}{\mathbb{T}}
\newcommand{\TT} { {\cal T }}
\newcommand{\FF} { {\cal F }}

\def\build#1_#2^#3{\mathrel{
\mathop{\kern 0pt#1}\limits_{#2}^{#3}}}
\newcommand{\UU}{{\cal U}}

\def\cq{$\hfill \square$}

\def\d{{\rm d}}

\def\e{{\cal E}}

\def\y{{\cal Y}}

\def\O{{\cal O}}

\def\eps{\epsilon}

\def \floor#1{\lfloor#1\rfloor}

\newcommand{\ssigma}{{\boldsymbol{\sigma}}}

\newcommand{\ms}{\mathbbm{s}}

\newcommand{\M}{\mathscr{M}}
\newcommand{\dgh}{{\rm d}_{\rm GH}}

\def\ov{\overline}

\def\la{\longrightarrow}

\def\proof{\noindent{\bf Proof. }}

\def \sint{\bigcup_{\ms \in \T}}

\renewcommand{\t}{\mathbbm{t}}
\newcommand{\ee}{\mathbbm{e}}
\newcommand{\ind}{\mathbbm{1}}

\newtheorem{thm}{Theorem}
\newtheorem{lmm}{Lemma}
\newtheorem{prp}{Proposition}
\newtheorem{defn}{Definition}
\newtheorem{crl}{Corollary}

\usepackage{fancyhdr}

\usepackage[outerbars]{changebar}

\begin{document}

\selectlanguage{english}

\begin{center}
\Large{\bf The CRT is the scaling limit of unordered binary trees}
\large ~\\
{\[\begin{array}{lcl}
 \textsf{Jean-Fran\c{c}ois Marckert}&~~&\textsf{Gr\'egory Miermont } \\
\textrm{CNRS, LaBRI}&~~ & \textrm{CNRS \& DMA}\\
\textrm{Universit\'e Bordeaux}&~~ &\textrm{Ecole Normale Sup\'erieure}\\
\textrm{351 cours de la Lib\'eration}&~~& \textrm{45 rue d'Ulm}\\
 \textrm{33405 Talence cedex}& ~~ &\textrm{F-75230 Paris Cedex 05} 
 \end{array}\]
}

\end{center}

\thispagestyle{empty}

\begin{abstract}
We prove that a uniform, rooted unordered binary tree with $n$
vertices has the Brownian continuum random tree as its scaling limit
for the Gromov-Hausdorff topology. The limit is thus, up to a constant
factor, the same as that of uniform plane trees or labeled trees. Our
analysis rests on a combinatorial and probabilistic study of
appropriate trimming procedures of trees.
\end{abstract}

\section{Introduction}\label{sec:introduction}

The Brownian Continuum Random Tree (CRT), introduced by Aldous
\cite{aldouscrt91}, is a natural object that arises in various
situations in Probability Theory. It is known to be the universal
scaling limit for conditioned critical Galton-Watson trees with finite
variance offspring distribution \cite{aldouscrt93,mm01,duq02}, or of
random labeled trees on $n$ vertices (Cayley trees)
\cite{aldouscrt91,jpmc97b,almpiep}.

Several distinct proofs for the convergence of discrete trees towards
the CRT exist in the literature, taking advantage of the specific
aspects of the considered models, which in passing yield various
equivalent constructions of the CRT. For instance, a
specific\footnote{This can also be obtained in the framework of
  Galton-Watson trees by choosing the particular Poisson offspring
  distribution} proof of the convergence of Cayley trees rests on the
fact that a uniform Cayley tree is a uniform spanning tree of the
complete graph, which can be constructed {\em via} the Aldous-Broder
algorithm and leads to a limiting ``stick-breaking construction'' of
the CRT. On the other hand, the convergence of conditioned
Galton-Watson trees to the CRT can be obtained by appropriate
encodings of trees by random walks. In \cite{aldouscrtov91}, Aldous
conjectures that many other models of trees, for which it is harder to
have a good probabilistic understanding, also have the CRT as a
continuum limit.

In this article, we will focus on one of these models, namely, the
family of rooted binary unordered trees, considered as graph-theoretic
trees without planar or labeled structure.  The main goal of this
paper is to prove that a uniformly chosen rooted binary unordered tree
with $n$ leaves converges, after renormalization of distances by
$\sqrt{n}$, to a constant multiple of the CRT.

\begin{figure}[htbp]
\centerline{\includegraphics[height=1.6cm]{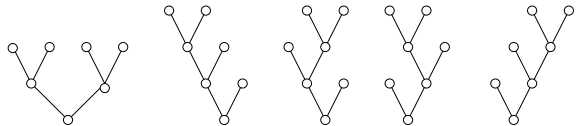}}
\caption{Representation of the 5 rooted binary plane trees with 4
  leaves. There are 2 (rooted) binary trees with 4 leaves (the 4 last
  ones on the picture representing the same unordered tree). }
\label{NPT}
\end{figure}

From a combinatorial perspective, models of plane trees are very close
to the case studied here: their generating functions' singularities
are of the same type (a square-root singularity), and heuristically,
combinatorial families with generating functions bearing the same
singularity type should have a similar continuous limit.
Nevertheless, the probabilistic methods developed for plane or labeled
trees are not valid anymore for unordered trees, and another method
has to be developed to derive their scaling limits. The present work
is a step in this direction, for the particular case of binary,
rooted, unordered trees.

We mention that in a very recent paper, another model of rooted
unlabeled, unordered trees has been investigated by Drmota \&
Gittenberger \cite{drmgit09}: this is the model of P\'olya trees where
there is no restriction on the authorized degrees. They compute the
limiting height and limiting profile of these trees (with $n$ nodes,
when $n\to +\infty$) and show that, up to constant factors, they
respectively converge to the maximum of the Brownian excursion, and
the local time of the Brownian excursion (which is similar to the case
of binary plane trees). Even if this result is of a different kind
from ours, and does not imply the convergence of the rescaled tree to
the CRT for the Gromov-Hausdorff topology, it gives an important
indication that it should be the case.

\subsection{Plane trees}\label{sec:plane-trees}

\begin{figure}[htbp]
\begin{center}
\psfrag{t}{$\t$}
\psfrag{v}{$\varnothing$}
\psfrag{1}{1}
\psfrag{11}{11}
\psfrag{111}{111}
\psfrag{112}{112}
\psfrag{121}{121}
\psfrag{122}{122}
\psfrag{12}{12}
\psfrag{21}{21}
\psfrag{22}{22}
\psfrag{2}{2}
\includegraphics[height=2.5cm]{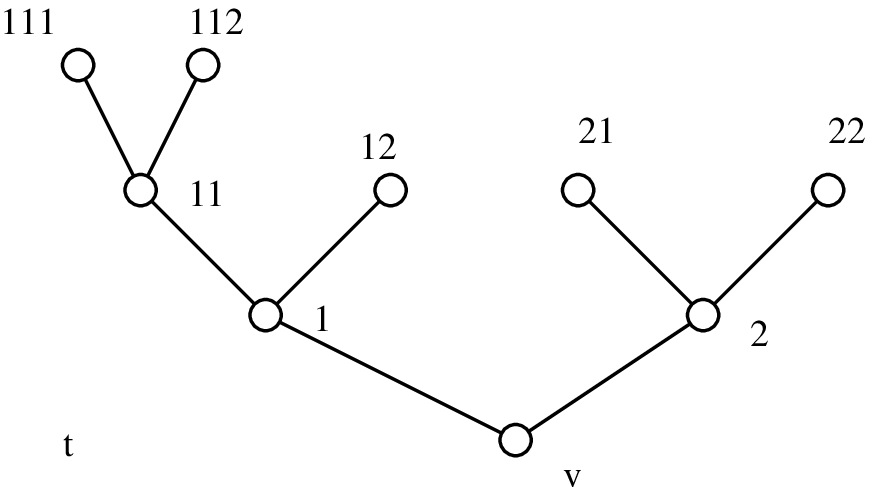}
\end{center}
\caption{A plane tree $\t$ such that $L(\t)=\{111,112,12,21,22\}$,
  $I_2(\t)=\{\varnothing\}$, $I_0(\t)=\{11,2\}$, $I_1(\t)=\{1\}$,
  $S(\t)=\{12\}$ .}
\label{PT2}
\end{figure}

\begin{figure}[htbp]
\begin{center}
\psfrag{t}{$\t$}
\psfrag{t'}{$\t'$}
\psfrag{v}{$\varnothing$}
\psfrag{1}{1}
\psfrag{11}{11}
\psfrag{111}{111}
\psfrag{112}{112}
\psfrag{121}{121}
\psfrag{122}{122}
\psfrag{12}{12}
\psfrag{21}{21}
\psfrag{22}{22}
\psfrag{2}{2}
\includegraphics[height=2.5cm]{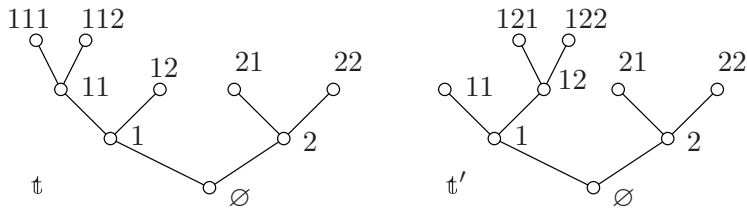}
\end{center}
\caption{Two plane trees $\t$ and $\t'$ that are equal as trees.}
\label{PT}
\end{figure}

Let us recall the standard definition of plane trees (see Figures
\ref{NPT}, \ref{PT2} and \ref{PT}).  Let $\UU=\bigsqcup_{n\geq 0}\N^n$
be the set of words with integer letters, where $\N=\{1,2,\ldots\}$
and $\N^0=\{\varnothing\}$.  For $u\in \UU$, let $|u|$ be its length,
and we generally note $u=u_1\ldots u_n$ with $n=|u|$. The
concatenation of the words $u$ and $v$ is written $uv$, and we write
$u\preceq v$ if $u$ is a prefix of $v$, meaning that there exists
$w\in\UU$ such that $uw=v$. This defines a partial order on $\UU$. If
$A\subset \UU$ and $u\in \UU$, we let
$$uA=\{uv:v\in A\}\, .$$

If $u,v\in\UU$, we let $u\wedge v$ be the longest common prefix to $u$
and $v$.  The set $\UU$ is endowed with the lexicographical order $\leq$: 
we have $u\leq v$ if $u\preceq v$ or if $u\wedge v$ is a strict
prefix of $u$ and $v$ such that $u_{|u\wedge v|+1}<v_{|u\wedge v|+1}$ 
(on $\mathbb{N}$). The order $(\UU,\leq)$ is total.

\begin{defn}
  A {\em rooted plane tree} is a finite subset $\t\subset \UU$
  containing $\varnothing$, such that if $ui\in \t$ with $u\in \UU$
  and $i\in \N$, then
\begin{itemize}
\item $uj\in \t$ for $1\leq j\leq i$, and
\item $u\in\t$.
\end{itemize}
(See Figure \ref{PT2}).
The elements of $\t$ are called vertices, and $\varnothing$ is called
the root. A vertex $ui\in\t$, with $u\in\UU$ and $i\in \N$, is called a child
of $u$. Their number is denoted by $c_u(\t)=\sup\{i\geq
1:ui\in\t\}\in \Z_+=\{0,1,2,\dots\}$. The notion of brothers, ancestors, descendants 
are induced by that of child as in the standard life.
The length $|u|$ is equally called {\em height}
of $u$.

In a plane tree $\t$, the subtree of $\t$ rooted at $u\in \t$ is the
plane tree $\t_u=\{v\in \UU:uv\in \t\}$.
\end{defn}

A plane tree has a representation as a plane graph, where each vertex
is linked by an edge to its children, which are ordered from left to
right in lexicographical order (as done on Figure \ref{PT2}).

We say that a rooted plane tree $\t$ is binary if $c_u(\t)\in\{0,2\}$
for every $u\in \t$. The vertices having no children are called the
{\em leaves} the set of which is denoted by $L(\t)$, while
$I(\t)=\t\setminus L(\t)$ denotes the set of {\em internal vertices}.
In turn, the internal vertices can be partitioned into three sets
$$I(\t)=I_0(\t)\sqcup I_1(\t)\sqcup I_2(\t)\, ,$$ where for
$i\in\{0,1,2\}$, $I_i(\t)$ is the set of internal vertices having $i$
children being themselves internal vertices (see Figure \ref{PT2}).
We let
$$S(\t)= L(\t)\setminus \{u1,u2:u\in I_0(\t)\}\, ,$$
the set of leaves that are children of the vertices of $I_1(\t)$. These
leaves will be called ``skeleton leaves'' of $\t$. 

 We let $|\t|$ denote the number of leaves of $\t$.  We let $\T$ be
 the set of binary rooted plane trees, and $\T_n=\{\t~:~
 |\t|=n\}\subset \T$ be the subset of those having $n$ leaves (hence
 $2n-1$ vertices and $2n-2$ edges).  It is a simple exercise to
 enumerate rooted binary plane trees via their generating
 functions. Namely, we have
$$\#\T_n=\frac{1}{n}\binom{2n-2}{n-1}\, ,$$ the $(n-1)$-th Catalan
 number, which also counts numerous families of combinatorial
 objects. Note the equivalent
\begin{equation}\label{eq:5}
\#\T_n\build\sim_{n\to\infty}^{}\frac{4^{n-1}}{\sqrt{\pi}n^{3/2}}\, .
\end{equation}

\subsection{Unordered trees}\label{sec:trees}

A rooted, binary unordered tree is a rooted, binary plane tree in
which the planar order has been ``forgotten'' (but not the root).
Two rooted binary plane trees $\t$ and $\t'$ are equivalent if $\t'$
can be obtained from $\t$ by exchanging the order of the various
children of $\t$, while preserving the rooted genealogical structure:
ancestors and brothers remain as such in the process.  Formally,
consider $\t\in \T$ and a family of permutations
$\ssigma=(\sigma_u,u\in I(\t))$ of the set $\{1,2\}$, indexed by the
internal vertices of $\t$. For $u=u_1\ldots u_n\in \t$, we let
$$\ssigma(u)=\sigma_\varnothing(u_1)\sigma_{u_1}(u_2)
\ldots\sigma_{u_1\ldots u_{n-1}}(u_n)\, ,$$ and
$\ssigma(\varnothing)=\varnothing$.
 Now for $\t'\in \T$ we say that
$\t\sim \t'$ if there exists such a $\ssigma$ with
$$\t'= \ssigma(\t)\, .$$ Note that there are at most $2^{|\t|-1}$
 elements in the class of $\t$, since $\#I(\t)=|\t|-1$.
\begin{defn}\label{sec:trees-1}
The quotient set $\bT=\T/\sim$ is called the set of (binary, rooted,
unordered) trees.  projection.  We also let $\bT_n=\T_n/\sim$ be those
trees which have $n$ leaves, and set $|\bt|=n$ for $\bt\in \bT_n$.
\end{defn}

It will be convenient to fix once and for all a section $s:\bT\to\T$
of the canonical projection $\T\to \bT$, so $s(\bt)$ is a choice among
all the planar representatives of $\bt$. In this way, there will be a
non-ambiguous way to consider a vertex of an unordered tree, the
number of its children, one of its edges, and so on.

The numbers $(\#\bT_n,n\geq 1)$ are the {\em Etherington-Wedderburn
  numbers}, the first ones are $1,1,1,2,3,6,11,23,46,98,207,451,...$
(referred to as A001190 in Sloane's On-Line Encyclopedia of Integer
Sequences).
No closed formula is available, but the asymptotic enumeration of
elements of $\bT$ has been studied by Otter \cite{otter48}, using the
properties of the generating function for trees\footnote{There is a
  slight difference in the article of Otter, as he considers trees
  with vertex-degrees either equal to $1$, $2$ or $3$. It is easy to
  see that this is actually equivalent to the family of binary
  unordered trees, by adding ``ghost'' edges at every vertex with only
  one child.}.  Let $\mathfrak{T}(z)=\sum_{n\geq 1}\#\bT_n\, z^n$ be
the generating function for rooted binary trees, counted according to
their number of leaves. A decomposition of a tree at its root into a
multiset of two subtrees yields the formula
\begin{equation}\label{eq:1}
\mathfrak{T}(z)=z+\frac{\mathfrak{T}(z^2)+\mathfrak{T}(z)^2}{2}\, .
\end{equation}
One shows that the radius of convergence $\varrho$ of $T$ satisfies
\begin{equation}\label{eq:12}
\mathfrak{T}(\varrho)=1\, .
\end{equation}
A detailed analysis of the recursive equation (\ref{eq:1})
entails
\begin{equation}\label{eq:3}
\#\bT_n\build\sim_{n\to\infty}^{}\frac{\bc}{2\sqrt{\pi}
  n^{3/2}\varrho^n}\, ,
\end{equation}
where 
\begin{equation}\label{eq:4}
\bc:=\sqrt{2\varrho+2\varrho^2 \mathfrak{T}'(\varrho^2)}=1.1300337\ldots\, ,
\end{equation}
see \cite[Note VII.22, p.477]{FlSe}.
The equivalent (\ref{eq:3}) is similar to (\ref{eq:5}), which hints at
the similarity of structure of large trees with large planar trees.

\medskip

\noindent{\bf Important notice. } From this section onwards, we will
only be interested in rooted, binary, plane or unordered trees, and
the words ``rooted, binary'', will {\em always} be implicit when
dealing with the two combinatorial families of trees we have just
defined. To avoid confusion, from now on, 
\begin{itemize}\item
a {\em plane tree} stands
for a rooted, binary plane tree, 
\item
while a {\em tree} is a rooted,
binary, unordered tree.
\end{itemize}

\subsection{Gromov-Hausdorff topology}\label{sec:grom-hausd-topol}

It is natural to consider trees (plane or not) as finite metric
spaces, by endowing the set of their vertices with the usual graph
distance. Formally, for $\t\in \T$ and $u,v\in \t$, we let
$$d_\t(u,v)=|u|+|v|-2|u\wedge v|\, .$$ In the sequel, we will often
improperly identify $\t$ with the metric space $(\t,d_\t)$, or even
with its isometry class, i.e.\ the collection of metric spaces that
are isometric to it.  It is straightforward that $(\t,d_\t)$ and
$(\t',d_{\t'})$ are isometric spaces whenever $\t\sim \t'$, so that a
tree $\bt$ will also denote the isometry class of the metric space
$(s(\bt),d_{s(\bt)})$. We will adopt the notation $a(X,d)=(X,ad)$ for
$(X,d)$ a metric space (identified with its isometry class) and
$a>0$. Therefore, the notation $a\t,a\bt$ will stand for the metric
spaces $\t,\bt$ with distances multiplied by $a$.

There is a common way to define a topology (even a metric) on spaces
of metric spaces. Such topologies have been developed in Geometry for
the last 30 years, following the ideas of Gromov
\cite{gromov99}. Their use in Probability in the context of random
real trees has been popularized by Evans and his coauthors
\cite{MR1774068,evpiwin,evanswinter}, and has been applied
successfully in various situations, for instance by Duquesne \& Le
Gall \cite{duqlegprep}, Le Gall \cite{legall06}.
We let $\M$ be the set of
isometry classes of compact metric spaces, which we endow with the
{\em Gromov-Hausdorff distance} defined by
\begin{equation}\label{eq:6}
\dgh((X,d),(X',d'))=\inf_{\phi,\phi'}\delta_H(\phi(X),\phi'(X'))\,
,
\end{equation}
the infimum being taken over the set of isometric embeddings
$\phi,\phi'$ from $(X,d),(X',d')$ into a common metric space
$(Z,\delta)$, and where $\delta_H$ is the usual Hausdorff distance
between compact subsets of $Z$. The reader can consult \cite[Chapter
  7]{burago01} for basic properties of this distance, which turns $\M$
into a Polish metric space (this is a straightforward extension of the
proof of \cite[Theorem 1]{evpiwin}). 

Of particular importance for our purposes is the subset
$\mathscr{T}\subset \mathscr{M}$ of $\R$-trees, i.e.\ of compact
spaces $(\TT,d)$ such that for every pair of points $x,y\in \TT$,
\begin{itemize}
\item 
there exists an isometry
  $\varphi_{x,y}:[0,d(x,y)]\to \TT$ with $\varphi_{x,y}(0)=x$ and
  $\varphi_{x,y}(d(x,y))=y$.
\item for every continuous injective $q:[0,1]\to \TT$ with $q(0)=x$
  and $q(1)=y$, it holds that $q([0,1])=\varphi_{x,y}([0,d(x,y)])$.
\end{itemize}
The first property says that $(\TT,d)$ is a {\em geodesic space}, and
the second is a ``tree property'' that there is a unique way to travel
between two points without backtracking. The set $\mathscr{T}$ is
closed in $\mathscr{M}$.

We should mention that when dealing with the convergence of plane
trees, a natural topology is inherited from the uniform distance
between (rescaled) contour processes, as defined later in Section
\ref{sec:case-discrete-trees}. This topology is stronger than the the
Gromov-Hausdorff topology, since contour processes encode also the
lexicographical order of the vertices. For unordered trees, the
Gromov-Hausdorff distance is very natural since it is by essence a
distance on metric spaces, which ignores any extra structure, like
planarity.

Finally, we stress that there is a variant of the Gromov-Hausdorff
distance on pointed metric spaces \cite{evpiwin}, which would be
natural in our context since we are interested in rooted trees. It is
straightforward to check that the results to come remain valid in this
setting, by taking the root vertex as the distinguished point.

\subsection{The Brownian continuum random tree}\label{sec:brown-cont-rand}

In order to state our main result, it remains to briefly describe the
CRT, which arises as the scaling limit of plane trees, as we now
recall.

\begin{figure}[htbp]
\psfrag{0}{0}\psfrag{x}{$x$}
\psfrag{1}{$\sigma_f$}
\psfrag{s}{$s$}
\psfrag{s'}{$s'$}
\psfrag{t}{$t$}
\psfrag{f(s)}{$f(s)$}
\psfrag{f(t)}{$f(t)$}
\psfrag{m_f}{$\check{f}(s,t)$}
\centerline{\includegraphics[height= 3.6 cm]{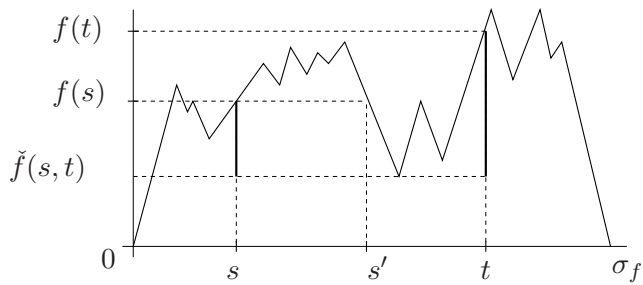}}
\caption{\label{abst-tree}Graph of a function $f$ from $\e$. 
In this example $s\equiv_f s'$ and the distance $d_f(s,t)=d_f(s',t)=f(s)+f(t)-2\,\check{f}(s,t)$.}
\end{figure}

Let ${\cal E}$ be the set of continuous functions $f$ defined on an
interval $[0,\sigma_f]$, that are non-negative, and satisfy
$f(0)=f(\sigma_f)=0$.  It is a complete space when endowed with the
distance
$${\rm D}(f,g)=\sup_{t\geq 0}|f(t\wedge \sigma_f)-g(t\wedge
\sigma_{g})|\, .$$ With every $f\in \e$, we can associate an
$\R$-tree, following \cite{duqlegprep}.  Define a function $d_f$ on
$[0,\sigma_f]^2$ by
$$d_f(s,t)=f_s+f_t-2\check{f}_{s,t}\, ,$$ where by definition
$$\check{f}_{s,t}=\inf_{s\wedge t\leq u\leq s\vee t}f_u\, .$$ It
is easy to see that $d_f$ is a pseudo-distance on $[0,\sigma_f]$. It
is not a distance as it does not separate points, so we let $s\equiv_f
t$ if $d_f(s,t)=0$, defining an equivalence relation on
$[0,\sigma_f]$ (see Figure \ref{abst-tree}).

The quotient space $\TT_f:=[0,\sigma_f]/\equiv_f$ is endowed with the
distance induced by $d_f$, which we still call $d_f$.  The canonical
projection $[0,\sigma_f]\to \TT_f$ being obviously continuous, the
target space is compact.  The (isometry class of the) space $\TT_f$ is
thus an element of $\mathscr{M}$, and turns out to be an element of
$\mathscr{T}$. The space $\TT_f$ is naturally {\em rooted},
i.e.\ comes with a distinguished point $\rho$, which is the
$\equiv_f$-equivalence class of $0$.  Finally, the mapping $f\mapsto
\TT_f$ is continuous from $(\e,{\rm D})$ to $(\mathscr{T},\dgh)$.

Now, let $(\ee_t,0\leq t\leq 1)$ denote the standard normalized
Brownian excursion (\cite[Chapter XII]{revyor}).  The CRT is the random isometry 
class $\TT_{2\ee}$,
defining a random variable in $(\mathscr{T},\dgh)$. Note the
convention to use {\bf twice} the distance $d_\ee$. This is only a
matter of convenience, and some references call ``CRT'' the metric
space $\TT_\ee$.

\subsection{Main result}\label{sec:brownian-crt}

We are now ready to state our main theorem. As a motivation, we start
with a well-known result.  Let $\P_n$ be the uniform distribution on
$\T_n$. The following result is a re-interpretation in terms of
$\R$-trees of a result of Aldous, see \cite{aldouscrt93} and
\cite{duqlegprep}.

\begin{prp}\label{sec:brownian-crt-2}
Let $\TT_n$ be a random variable with distribution $\P_n$. Then
$$\frac{1}{\sqrt{2n}}\TT_n\build\la_{n\to\infty}^{(d)}\TT_{2\ee}\, ,$$
where the convergence in distribution holds with respect to the
Gromov-Hausdorff topology on $\mathscr{M}$. 
\end{prp}

Our main result is an analog of the previous statement for trees. We
let $\bP_n$ be the uniform distribution on $\bT_n$. 

\begin{thm}
\label{sec:bc1}
Let $\TT_n$ be a random variable with distribution $\bP_n$. Then
$$\frac{\bc}{\sqrt{2n}}\TT_n\build\la_{n\to\infty}^{(d)}\TT_{2\ee}\,
,$$ where the convergence in distribution holds with respect to the
Gromov-Hausdorff topology on $\mathscr{M}$, and where $\bc$ is given
in (\ref{eq:4}).
\end{thm}

In particular, one will notice that a tree with law $\bP_n$ is
typically about $\bc$ times smaller (indeed $\bc$ is greater than 1)
than a typical $\P_n$-distributed plane tree (in passing, a quick
computation based on Figure \ref{NPT} shows that the mean height of a
random variable with distribution $\P_4$ is 14/5, while it equals
$5/2$ for a random variable with distribution $\bP_4$). \par

\subsection{Contents of the paper and strategy of the proof of Theorem \ref{sec:bc1}}

In Section \ref{sec:case-plane-trees} we present a trimming operation
for plane and unordered trees: the idea is to keep in a large plane
tree with $n$ leaves, only the vertices having at least $\floor{`e n}$
descendants, which in turn can be encoded as a plane tree with
edge-lengths, called the $\floor{\eps n}$-skeleton.  A similar
operation is defined on (unordered) trees in Section
\ref{sec:case-trees}.

This operation is the main tool in our study. Asymptotically, the
$\floor{\eps n}$-skeletons of $\P_n$ and $\bP_n$-distributed random
variables will turn out to be directly comparable thanks to two very
similar ``local limit theorems'' (Propositions
\ref{sec:consequences-2} and \ref{sec:proof-theor-refs-1}), which give
the scaling limits for the densities of these $\lf \eps
n\rf$-skeletons. These are obtained via combinatorial arguments
developed in Sections \ref{sec:case-plane-trees} and
\ref{sec:case-trees}.

These local limit theorems form the cornerstone of our study, and they
are, at the intuitive level, the main explanation of the similar
asymptotic behaviour of both families of trees. However, these results
are not sufficient by themselves to entail Theorem
\ref{sec:bc1}. First, because the pointwise convergence of the density
of a random variable is not sufficient to get convergence in
distribution. We need to check that the limiting formulas in
Propositions \ref{sec:consequences-2} and \ref{sec:proof-theor-refs-1}
define indeed probability distributions, i.e.\ to ensure that no mass
disappears, or goes into a singular part. A second problem is that the
convergence of the $\floor{`e n}$ skeleton does not imply immediately
the convergence of the non-trimmed trees under $\bP_n$: we have to
rule out the possibility that the parts of the trees that have been
removed in the trimming procedure are thin ``hair'' of very large
diameter although they contain a small amount of leaves (at most
$\lf\eps n\rf$).
The identification of the limit is then a consequence of Proposition
\ref{sec:brownian-crt-2}.

Let us describe more precisely these steps. The local limit result,
Proposition \ref{sec:consequences-2}, concerning the $\floor{`e
  n}$-skeleton under $\P_n$, is shown to imply a convergence in
distribution in Section \ref{sec:density-psitt-its}. To this end, we
make a strong use of the convergence of the contour process of a
$\P_n$-distributed plane tree toward the Brownian excursion
(Proposition \ref{sec:trimm-read-cont-2}). This involves a careful
translation of our trimming operations in terms of operations on
contour processes, and more generally on excursion functions seen as
encoding $\mathbb{R}$-trees as explained in Section
\ref{sec:brown-cont-rand}. This is done in Sections
\ref{sec:trimm-read-cont} and \ref{sec:encoding-tt_fa-as}.

As a matter of fact, in the plane case, the convergence of contour
processes under $\P_n$ is sufficiently robust to entail that of the
$\floor{`e n}$-skeleton to an $`e$-trimmed version of the CRT,
described in term of the Brownian excursion (Corollary
\ref{sec:cont-funct-trimm-1} in Section \ref{sec:cont-prop-a}, and
Sections \ref{sec:scaling-limits} and \ref{sec:interpr-terms-psims}).
This gives the wanted convergence in distribution for $\floor{\eps
  n}$-trimming of a $\P_n$-distributed plane tree, and a simple
comparison argument (thanks to Propositions \ref{sec:consequences-2}
and \ref{sec:proof-theor-refs-1}), allows to prove that $\floor{\eps
  n}$-skeletons under $\bP_n$ converges also in distribution toward
the same limit, up to a constant factor (Section
\ref{sec:conv-a-skel}).

Finally, we show in Section \ref{sec:proof-theor-refs} that the
trimmed versions are not too far from the whole tree with high
probability: this part, which amounts to controlling the maximal
height of all the subtrees appearing in the difference between a tree
and its $\floor{`e n}$-skeleton relies on moment estimates for the
height of a $\bP_n$-distributed tree, obtained by Broutin \& Flajolet
\cite{BrFl08}.

\bigskip

\noindent{\bf Acknowledgments. }We are grateful to Mathilde Weill,
whose help and insights in the early stage of this work have been very
important. Grateful thanks are due to Nicolas Broutin and Philippe
Flajolet for the interest they have taken in the problem of estimating
the total height of a random $\bP_n$-distributed tree, and for keeping
us informed of the progress of their work, an extended abstract of
which can be found in \cite{BrFl08}. These results are indeed crucial
in obtaining a tightness argument of the kind of Lemma
\ref{sec:proof-theor-refs-2}.

\section{Trimming trees}\label{sec:trimming-trees}

Our main tool to prove Theorem \ref{sec:bc1} will be to use 
a mass-trimming of the tree from the leaves. The underlying idea is 
that the combinatorics of the resulting
tree are a lot more tractable, and show quite easily the universal
aspects of different tree models, at least for binary trees. In this
perspective, we are going to compare the trimming of plane trees and
trees with respective distributions $\P_n,\bP_n$.

\subsection{The case of plane trees}\label{sec:case-plane-trees}

\subsubsection{$a$-trimmed tree}
\label{tri-tree}
\begin{figure}[htbp]
\psfrag{4}{4}
\psfrag{0}{0}
\psfrag{1}{1}
\psfrag{5}{5}
\begin{center}
\includegraphics[height=3.6cm]{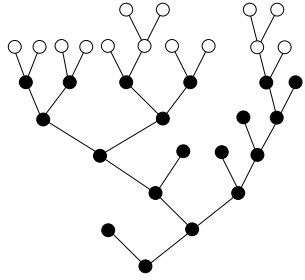}
\end{center}
\caption{A plane tree $\t$; the tree $\t[a]$ with $a=3$, i.e.\ those
  vertices whose fathers have 4 leaves at least among their
  descendants, is represented by the black vertices}
\label{tri}
\end{figure}

Let $\t$ be a plane
tree, and for $a\geq 0$, the {\em $a$-trimmed tree} is defined to be
$$\t[a]=\{u\in \t:|\t_v|> a\quad \forall\, v\prec u\}\, .$$ 
Of course, to check that $u\in \t[a]$ it suffices to verify that
$|\t_v|>a$ where $v$ is the parent of $u$. However, with our
definition, $\t[a]$ always contains $\varnothing$, and obviously
defines a plane tree. For instance, $\t[0]=\t[1]=\t$ (see also Figure \ref{tri}). From now on we
assume that $a\geq 2$.

Let us discuss the reconstruction of a plane tree $\t$, starting from
its $a$-trimming $\t[a]$. By definition of $\t[a]$ it must hold that
$|\t_u|\leq a$ for every $u\in L(\t[a])$ (otherwise the children of $u$
would belong to $\t[a]$ as well), and moreover, if $u$ is such that
$u1,u2$ are both in $L(\t[a])$, then it must hold that
$|\t_{u1}|+|\t_{u2}|>a$ (otherwise $|\t_u|\leq a$ so $u1$ could not be
in $\t[a]$). The following lemma is a (straightforward) kind of
converse to this observation.
\begin{lmm}
  Let $\t_0\in \T$, and $(\t_{(u)},u\in L(\t_0))\in \T^{L(\t_0)}$ be a family 
  of plane trees indexed by the leaves of $\t_0$, such that 
  \begin{itemize}
  \item $|\t_{(u)}|\leq a$ for every $u\in \t$,
  \item $|\t_{(u1)}|+|\t_{(u2)}|>a$ whenever $u1,u2\in L(\t)$, 
  \end{itemize}
therefore $\t_0=\t[a]$, where 
\begin{equation}\label{eq:19}
\t=\t_0\cup \bigcup_{u\in L(\t_0)}u\t_{(u)}\, .  
\end{equation}
\end{lmm}

Consequently, a plane tree can be recovered in a one-to-one fashion
from its $a$-trimming $\t[a]$ and an appropriate family of plane trees
indexed by the leaves of $\t[a]$. Define
$$\T_{\leq a}=\{\t\in \T:|\t|\leq a\},$$ the set of trees having at most $a$ leaves, and 
$$\F_{a}=\{\ff=(\t,\t'): \t,\t'\in \T,|\t|\vee |\t'|\leq
a,|\t|+|\t'|>a\}\, .$$ A pair of the form $\ff=(\t,\t')$ is called a
forest with two tree components, namely $\ff_1=\t,\ff_2=\t'$ with
natural notations. The size $|\ff|$ is then defined as
$|\ff|=|\t|+|\t'|$, and by convention, for $u\in {\cal U}$ and
$\ff=(\t,\t')$, we let $u\ff$ be the set $\{u\}\cup u1\t\cup u2\t'$.

With this notation, we can rewrite (\ref{eq:19}) as
$$\t=\t_0\cup \bigcup_{u\in S(\t_0)}u\t_{(u)}\cup \bigcup_{u\in
  I_0(\t_0)}u\ff_{(u)}\, ,$$ where the $\t_{(u)}$ are in $\T_{\leq a}$
and the $\ff_{(u)}$ are in $\F_a$. In the forthcoming Lemma
\ref{sec:trimming-trees-1}, we will make one further adjustment by
ordering the elements of $S(\t_0)$ and $I_0(\t_0)$ in lexicographical
order, allowing us to label the families $\t_{(u)}$ and $\ff_{(u)}$ by
integers rather than elements of $\t$.\medskip

In order to study the shape of $\t[a]$ we introduce two operations on
plane trees. Their aim is to encode in some sense $\t[a]$: when
$a$ is large, the tree $\t[a]$ will have few true branching points
($\#I_2(\t[a])$ is small), and a lot of skeleton leaves. The
operations defined below in the context of general binary plane trees
will then be applied on $\t[a]$ in Section \ref{sec:consequences}.

\subsubsection{Contraction of a plane tree}
\label{tro}
We first define an operation of concatenation.
Assume given two plane trees $\t,\t'$, together with
labels $\bx=(x_u,u\in \t),\bx'=(x'_u,u\in \t')$ taking their values in some
set $E$ (typically, $E$ will be the set $\Z_+$). For $x\in E$, we
define the {\em concatenation} 
$$\Xi(x,(\t,\bx),(\t',\bx'))=(\t'',\bx'')\, ,$$ to be the tree having
its root $\varnothing$ marked by $x$, and having the two marked tree
$\t$ and $\t'$ as subtrees rooted at the children $1$ and $2$ of
$\varnothing$. Formally
$$\t''=\{\varnothing\}\cup 1\t\cup 2\t'\, \quad \in \T\, ,$$
and $\bx''=(x_u'',u\in \t'')$ is defined by
$$x''_\varnothing=x\, ,\quad x''_{1u}=x_u\, ,\ u\in \t\, ,\quad 
x''_{2u}=x'_u\, ,\ u\in \t'\,. $$

\begin{figure}[htbp]
\psfrag{4}{4}
\psfrag{0}{0}
\psfrag{1}{1}\psfrag{2}{2}
\psfrag{5}{5}\begin{center}
\includegraphics[height=3.2cm]{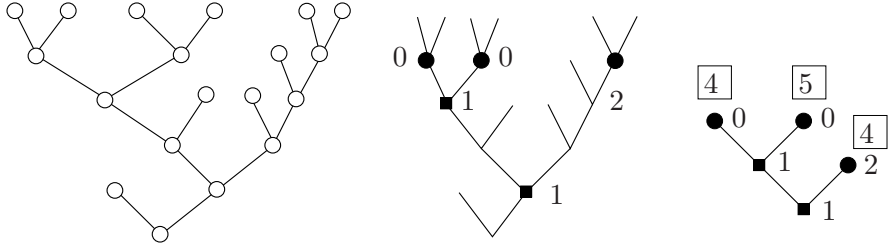}
\end{center}
\caption{The plane tree $\t'=\t[3]$ (where $\t$ is given in Figure \ref{tri}). 
The elements of $I_2(\t')$ (represented with squares), and elements of $I_0(\t')$ 
(represented by small black discs) are the vertices giving the structure of the projection 
$\Pi(\t')=(\ms,\bx)$, represented on the third picture.}
\label{fig:pi}
\end{figure}

Let us introduce the contraction application $\Pi$.  The action of
$\Pi$ amounts to suppressing all the leaves of $\t$, and then to
contracting the maximal chains of vertices with degree $2$ in the
resulting graph (while keeping its planar structure). This provides a
new labeled binary tree $\t'$, the label of a vertex being simply the
size of the chain that has been contracted underneath it (see Figure
\ref{fig:pi}).\par Let us be formal: for every $\t\in
\T\setminus\{\{\varnothing\}\}$ (so that $|\t|\geq 2$), we define the
                              {\em contraction} of $\t$ as
                              follows. Let $U(\t)$ be the vertex with
                              minimal height having degree 2 in $\t$
                              satisfying
$$\forall v\prec u\, , \quad |\t_{v1}|=1\ {\rm xor}\
|\t_{v2}|=1\, ,$$ where xor denotes exclusive or. It must be
understood that if no such $v$ as in the definition exists, then
$U(\t)=\varnothing$. Then the contraction of $\t\in \T$ is defined in
a recursive way as a plane tree with integer labels
$\Pi(\t)=(\ms,\bx)$ as follows:
\begin{itemize}
\item
If $|\t_{U(\t)1}|=|\t_{U(t)2}|=1$ then 
$\ms=\{\varnothing\}$ and $x_\varnothing=|U(\t)|$\,,
\item
otherwise, it must hold that $|\t_{U(\t)1}|\wedge |\t_{U(\t)2}|\geq
2$, and we set
\begin{equation}\label{Pit}
\Pi(\t)=\Xi(|U(\t)|,\Pi(\t_{U(\t)1}),\Pi(\t_{U(\t)2}))\,.
\end{equation}
\end{itemize}

This branching contraction procedure, together with the interpretation
given at the beginning of Section \ref{tro}, entails that the vertices
of $\Pi(\t)$ are naturally associated with a subset of the vertices of
$\t$ via an application $\pi_\t$. As one may observe on Figure
\ref{fig:pi}, the branching structure of $\Pi(\t)$ is given by the
relative positions of the vertices in $I_2(\t)\cup I_0(\t)$. More
precisely, we let $\pi_\t(\varnothing)=U(\t)$, and recursively, for
$u\in I(\Pi(\t))$, we set
\begin{equation}\label{eq:18}
\pi_\t(u1)=U(\t_{\pi_\t(u)1})\, ,\qquad
\pi_\t(u2)=U(\t_{\pi_\t(u)2})\, .
\end{equation}
It is easy to check that $\pi_\t$ is injective, with image $I_2(\t)\cup
I_0(\t)$.  More precisely, if $\ms=\Pi(\t)$, then $\pi_\t(I(\ms))=I_2(\t)$
and $\pi_\t(L(\ms))=I_0(\t)$. Moreover, $\pi_\t$ is increasing for the
lexicographical order.

Noticing that the plane trees contracted to the tree $(\{\varnothing\},x)$
are exactly those of the form $\{v,v1,v2:v\preceq u\}$ for a word
$u\in \{1,2\}^x$, and using an induction argument, one may prove straightforwardly the following Lemma:
\begin{lmm}\label{sec:trimming-trees-2}
  Let $\ms$ be a binary plane tree and $\bx=(x_u,u\in \ms)$. Then
$$\#\Pi^{-1}(\ms,\bx)=2^{|\bx|}\, ,$$
where $|\bx|=\sum_{u\in \ms}x_u$.  
\end{lmm}

\subsubsection{Encoding of the $a$-trimming of $\t$: the $a$-skeleton of $\t$}

We now go one step further in our simplification of trees.  Let
$(\ms,\bx)=\Pi(\t[a])$ be the contracted tree associated with the
$a$-trimming of $\t$. It is implicit that $\Pi(\t[a])$ is indeed
well-defined, meaning that $|\t[a]|\geq 2$, and this is equivalent to
the fact that $|\t|>a$, which is thus assumed.

For our purposes, it will be useful to assign
to the leaves of $\ms$ a {\em second integer label}, giving the total
number of descendant leaves they had in the initial tree (corresponding
to the size of the corresponding forest of pairs of trees).  Formally, for
$u\in L(\ms)$, 
$$y_{u}=|\t_{\pi_{\t[a]}(u)}|\, ,$$ 
where the mapping $\pi_{\t[a]}$ is as defined in (\ref{eq:18}),
and note that if $u\in L(\ms)$, it must hold that $\pi_{\t[a]}(u)$
is an element of $I_0(\t[a])$.  Also, it holds that $y_v>a$ for every
$v\in L(\ms)$. 

\begin{defn}\label{sec:case-plane-trees-1}
The triple $(\ms,\bx,\by)$ is called the
$a$-skeleton of $\t$, and is denoted by ${\rm Sk}_a(\t)$.
\end{defn}

The general $a$-{\em skeleton} is thus a triple $(\ms,\bx,\by)$ such
that
\begin{itemize}
\item
$\ms\in \T$,
\item
$\bx=(x_u,u\in \ms)$ is a vector of non-negative integers,
\item
$\by=(y_u,u\in L(\ms))$ is a vector of integers in
  $(a,2a]$.
\end{itemize}

Summing up our study so far, we obtain the following result.

\begin{lmm}\label{sec:trimming-trees-1}
  Any (binary) plane tree $\t$ such that $|\t|>a$ can be obtained
  from the following data, in a one-to-one way:
\begin{enumerate}
\item an $a$-skeleton $(\ms,\bx,\by)$
\item one of the $2^{|\bx|}$ elements $\t_0\in \Pi^{-1}(\ms,\bx)$,
 \item two sequences $(\t(i),1\leq i\leq |\bx|)\in \T_{\leq
    a}^{|\bx|}$ and $(\ff(i),1\leq i\leq |\ms| )\in \F_a^{|\ms|}$ such
  that
$$|\ff(j)|=y_{v_j}\, ,\qquad 1\leq j\leq |\ms|\, ,$$
where $v_1\leq \ldots\leq v_{|\ms|}$ is the ordered list of leaves of
$|\ms|$
\end{enumerate}
by the formula:
$$\t=\t_0\cup \bigcup_{i=1}^{|\bx|} u_i\t(i)\cup
\bigcup_{j=1}^{|\ms|}l_j\ff(j)\, ,$$ where $u_1\leq \ldots \leq
u_{|\bx|}$ and $l_1\leq \ldots\leq l_{|\ms|}$ respectively denote the
ordered lists of elements of $S(\t_0)$ and $I_0(\t_0)$, so in
particular, $l_i=\pi_{\t_0}(v_j),1\leq j\leq |\ms|$.
\end{lmm}

\subsection{The case of trees}\label{sec:case-trees}

The procedure we adopt to trim a tree is similar to the one for plane
trees. The reconstruction properties of a tree from its trimming are
not as smooth as in the plane case though, as we will now see. First,
we focus on what we call a labeled tree. A labeled plane tree is a
pair $(\t,\bg)$ where $\t$ is a plane tree and $\bg:\t\to E$ is a
function from $\t$ to some set $E$.  We define a new equivalence
relation $\approx$ for labeled plane trees $(\t,\bg)$, by
$(\t,\bg)\approx (\t',\bg')$ if there exist permutations
$\ssigma=(\sigma_u,u\in \t)$ as in Sect.\ \ref{sec:trees} such that
$$\t'=\ssigma(\t)\, , \qquad g'_{\ssigma(u)}=g_{u}\, ,\qquad u\in
\t\, .$$ 

A labeled tree is just an equivalence class for $\approx$.  Such an
equivalence class has at most $2^{\#I(\bt)}=2^{|\t|-1}$
representatives. When this upper-bound is attained, we say that the
associated (rooted) labeled tree (or any of its plane representatives)
has {\em no symmetries}.  A simple criterion for having no symmetry is
given by the following intuitive result. We say that a labeled plane tree
$(\t,\bg)$ (or its $\approx$-equivalence class) is {\em good} if the
components of $\bg$ are pairwise distinct.

\begin{lmm}\label{sec:case-trees-1}
A good labeled plane tree $(\t,\bg)$ has no symmetries.
\end{lmm}

\proof Endow the label set $E$ with any particular linear order. With
any good labeled tree, one associates a particular planar
representative $(\t,\bg)$ in which for every pair of vertices born
from the same parent, the one with larger label is always the first
child. The other planar representatives are obtained by deciding the
order in which the two children of an interior vertex should be
displayed, knowing that they can always be distinguished, as they have
different labels, so there are exactly $2^{|\bt|-1}$ of them.  \cq 

\medskip

Next, the reader will easily be convinced that
$\Pi$ and $\t\mapsto \t[a]$ are class operations, meaning that if
$\t\sim \t'$ then
$\t[a]\sim \t'[a]$, and 
$\Pi(\t)\approx \Pi(\t')$. 
This validates the following
\begin{defn}
   We say that a plane tree $\t$ with $|\t|>a$ is
{\em $a$-good} if it satisfies the following properties
\begin{itemize}
\item the labeled tree $\Pi(\t[a])$ is good, and
\item for every $u\in I_0(\t[a])$, it holds that $|\t_{u1}|\neq
  |\t_{u2}|$.
\end{itemize}
We say that a tree $\bt$ is $a$-good if one (hence any) of its planar
representatives is. 
\end{defn}

If $\bt$ is an $a$-good tree, we can choose in some non-ambiguous way
a plane $\approx$-representative of its skeleton. To be specific, we
will set
$${\rm Sk}_a(\bt):=(\ms,\bx,\by)
$$ to be the representative such that $x_{u1}>x_{u2}$ for every $u\in I(\ms)$, and
call it the $a$-skeleton of $\bt$. The general $a$-skeleton of an
$a$-good tree (we naturally call such objects $a$-good skeletons) is
thus a triple $(\ms,\bx,\by)$ such that
\begin{itemize}
\item
$\ms$ is a plane tree,
\item
$\bx=(x_u,u\in \ms)$ is such that $(\ms,\bx)$
  is good and $x_{u1}>x_{u2}$ whenever $u\in I(\ms)$, and
\item $\by=(y_u,u\in L(\ms))$ is such that $y_u\in (a,2a]$ for $u\in
  L(\ms)$.
\end{itemize}

The key observation is the following fact. 

\begin{lmm}\label{sec:case-trees-7}
  Let $\t\sim\t'$ be $a$-good plane trees such that $\t[a]=\t'[a]$,
  and such that for every $u\in I_0(\t[a])$, one has 
$$|\t_{u1}|>|\t_{u2}|\qquad \mbox{ and }\qquad 
  |\t'_{u1}|>|\t'_{u2}|\, .$$ Then for every $u\in L(\t[a])$, it holds
  that $\t_u\sim \t'_u$.
\end{lmm}

\proof Let $\ssigma=(\sigma_u,u\in \t)$ be such that
$\ssigma(\t)=\t'$. Assume that there exists $u\in I(\t[a])$ such that
$\sigma_u\neq 1$ is not the identity. Let $u$ be such a vertex of
$I(\t[a])$, with minimal height. Then $\sigma_v=1$ for all $v\prec u$,
so since $\ssigma(\t)=\t'$, it must hold that $\t_u\sim \t'_u$.

First assume $u\in I_0(\t[a])$, then by assumption, it holds that
$|\t_{u1}|>|\t_{u2}|$ and $|\t'_{u1}|>|\t'_{u2}|$. Since $\t_u\sim
\t'_u$ and $\sigma_u\neq 1$, this entails that $\t_{u1}\sim \t'_{u2}$
(and $\t_{u2}\sim \t'_{u1}$), a contradiction with the order on the
sizes of these trees.

A second possibility is that $u\in I_1(\t)$, so that $u1$ or $u2$ is
in $S(\t[a])$. Assuming e.g.\ that the first case holds, this says
that $u2$ is in $I(\t[a])$, but since $\sigma_u\neq 1$ this says that
$u1\in I(\t'[a])=I(\t[a])$ which is disjoint from $S(\t[a])$, a
contradiction.

The last possibility is that $u\in I_2(\t[a])$, and we let
$v=\pi_{\t[a]}(u)\in I(\ms)$, where $(\ms,\bx)=\Pi(\t[a])$. But the
fact that $\sigma_u\neq 1$ has the effect of switching the values of
$x_{v1},x_{v2}$ in the (good) labeled plane tree $(\ms,\bx)$. Since
these values are distinct, this contradicts $\t[a]=\t'[a]$.

This shows that $\sigma_u=1$ for every $u\in I(\t[a])=I(\t'[a])$. The
conclusion follows easily by noticing that $\t_u'=\ssigma^{(u)}(\t_u)$
for every $u\in I(\t[a])$, where by definition,
$\ssigma^{(u)}=(\sigma_{uv}:v\in I(\t_u))$.  \cq

\medskip

We now choose, once and for all, a section $\Pi_0^{-1}$ of $\Pi$,
i.e.\ a mapping from the set of $a$-skeletons to $\T$ such that
$\Pi\circ \Pi_0^{-1}$ is the identity.

Starting from an $a$-good tree $\bt$ with $a$-skeleton  
$(\ms,\bx,\by)$, we can always choose a plane representative $\t$ of
$\bt$ such that $\t[a]=\Pi_0^{-1}(\ms,\bx)$. There might be many such
choices of $\t$, but by Lemma \ref{sec:case-trees-7}, we know there is
a well-defined family of trees $\bt_{(u)},u\in L(\t_0)$, such that
$|\bt_{(u1)}|>|\bt_{(u2)}|$ whenever $u\in I_0(\t_0)$, and such that
$\bt$ is the $\sim$-equivalence class of the plane tree (recall that
$s:\bT\to \T$ is a section)
$$\t[a]\cup \bigcup_{u\in L(\t_0)}us(\bt_{(u)})\, .$$
More precisely, letting
$$\bT_{\leq a}=\{\bt\in \bT:|\bt|\leq a\}$$ and 
$$\bF_{a}=\{(\bt,\bt'): a\geq |\bt|> |\bt'|,\, |\bbf|:=|\bt|+|\bt'|>a\}\, ,$$ 
the trees $\bt_{(u)},u\in S(\t[a])$ lie in $\bT_{\leq a}$, while the
pairs $(\bt_{(u1)},\bt_{(u2)})$ for $u\in I_0(\t[a])$ lie in $\bF_a$ and
satisfy $|\bt_{(u1)}|+|\bt_{(u2)}|=y_{\pi_{\t[a]}(u)}$.  By
convention, for $\bbf=(\bt,\bt')\in \bF_a$ we let
$s(\bbf)=(s(\bt),s(\bt'))\in \F_a$.  Summing up our study, we get:

\begin{lmm}\label{sec:case-trees-2}
Any $a$-good tree $\bt$ can be obtained from the following data in a
one-to-one way:
\begin{enumerate}
\item a good $a$-skeleton $(\ms,\bx,\by)$ such that $x_{u1}>x_{u2}$
  for every $u\in I(\ms)$, 
\item
two sequences $(\bt(1),\ldots,\bt(|\bx|))\in \bT_{\leq a}^{|\bx|}$ and
$(\bbf(1),\ldots,\bbf(|\ms|))\in \bF_a^{|\ms|}$ such that
$$\bbf(j)=(\bt(j),\bt'(j))\quad \mbox{
  and }\quad |\bbf(j)|=y_{v_j}\, ,\qquad 1\leq j\leq |\ms|\, ,$$ where
$v_1,\ldots,v_{|\ms|}$ is the ordered list of elements of $L(\ms)$,
\end{enumerate}
by letting $\bt$ be the $\sim$-equivalence class of the planar tree
$$\t=\t_0\cup \bigcup_{i=1}^{|\bx|} u_is(\bt(i))\cup
\bigcup_{j=1}^{|\ms|}l_js(\bbf(j))\, .$$ where
$\t_0=\Pi^{-1}_0(\ms,\bx)$, and $u_1,\ldots,u_{|\bx|}$ resp.\
$l_1,\ldots,l_{|\ms|}$ are the ordered lists of elements of $S(\t_0)$
and $I_0(\t_0)$.
\end{lmm}

\section{Local limit theorems for $a$-skeletons of 
trees}\label{sec:consequences}

Let us draw the probabilistic consequences of our analysis of trimmed
trees, starting with plane trees.

\subsection{Plane trees}\label{sec:plane-trees-1}

Let $\mu=(\mu_k)_{k\geq 1}$ be the law of the size of
a Galton-Watson tree with offspring distribution
$(\delta_0+\delta_2)/2$; one has immediately
\begin{equation}\label{etet}
\mu_n=\frac{\#\T_n}{2^{2n-1}}\, ,\qquad n\geq 1\,, 
\end{equation}
and by Stirling's formula, the following equivalent holds:
\begin{equation}\label{eq:8}
\mu_n\build\sim_{n\to\infty}^{}\frac{1}{2\sqrt{\pi}n^{3/2}}\, .
\end{equation}
We recall that conditioned on having $n$ vertices, a Galton-Watson
tree with offspring distribution $(\delta_0+\delta_2)/2$ is uniform on
$\T_n$. 
Consider a probability space
$(\Omega,\FF,P)$ on which is defined an i.i.d.\ sequence
$X_1,X_2\ldots$ of $\mu$-distributed random variables.

\begin{lmm}\label{sec:consequences-1}
  Let $(\ms,\bx,\by)$ be an $a$-skeleton. Then, seeing
  ${\rm Sk}_a:\t\mapsto {\rm Sk}_a(\t)$ as a random variable,
\begin{eqnarray}
\lefteqn{\P_n({\rm
    Sk}_a=(\ms,\bx,\by))}\\ &=&\frac{1}{2^{2|\ms|-1}\mu_n}
P\left(\sum_{i=1}^{|\bx|}X_i=n-|\by|,\, \max_{1\leq i\leq |\bx|}
X_i\leq a\right)\prod_{u\in L(\ms)}P(X_1+X_2=y_u,X_1\vee X_2\leq a).
\nonumber
\end{eqnarray}
\end{lmm}

\proof By Lemma \ref{sec:trimming-trees-1}, the number of plane trees
admitting $(\ms,\bx,\by)$ for $a$-skeleton equals
$$2^{|\bx|}\#\left\{(\t(1),\ldots,\t(|\bx|))\in \T_{\leq
    a}^{|\bx|}:\sum_{i=1}^{|\bx|}|\t(i)|=n-|\by|\right\} \prod_{u\in
  L(\ms)}\#\{\ff\in \F_{a}:|\ff|=y_u\}\, ,$$ which
equals by (\ref{etet})
\begin{eqnarray*}
\lefteqn{2^{|\bx|+2n-2|\by|-|\bx|}\prod_{j=1}^{|\ms|}2^{2y_i-2}}\\
&\times &P\left(\sum_{i=1}^{|\bx|}X_i=n-|\by|\, ,\, \max_{1\leq i\leq |\bx|}
X_i\leq a\right)\prod_{u\in L(\ms)}P(X_1+X_2=y_u,X_1\vee X_2\leq a)\, ,
\end{eqnarray*}
giving the result after dividing by $\#\T_n$.  \cq

\medskip

From this, we deduce an important limiting result (Proposition
\ref{sec:consequences-2}), which can be thought of a kind of ``local
limit theorem'' for the $\lf \eps n\rf$ skeleton of a
$\P_n$-distributed element. 

From now on, we fix a number $\eps\in (0,1)$ and will usually omit its
mention to allow lighter notations. The number $\eps$ will be allowed
to move in the further Section \ref{sec:proof-theor-refs} (also at
the end of Section \ref{sec:interpr-terms-psims}), but we will take
care that there is no ambiguity at this moment.

For $\ms\in \T$ and $M>1/\eps$, set
\begin{equation}\label{eq:13}
  I_M(\ms):=\left\{(\bx,\by)\in \R_+^{\ms}\times \R_+^{L(\ms)}:
    \left.\begin{array}{ll} x_u\in (M^{-1},M) & u\in \ms\\y_u\in
        (\eps+M^{-1},2\eps) & u\in L(\ms)
\end{array}\right.
\right\}\, .
\end{equation}
We also define the family of functions $a_k$ (for $k\geq 0$) on $(0,\infty)$ by letting $a_0\equiv 1$,
$a_1(x)=(2\sqrt{\pi}x^{3/2})^{-1}\ind_{x\geq \eps}$ and for $k\geq 2$,
\begin{equation}\label{eq:2}
a_k(x)=\frac{1}{(4\pi)^{k/2}}\int_{u_1,\ldots,u_{k-1}\geq
  \eps}\frac{\ind_{\{u_1+\ldots+u_{k-1}\leq x-\eps\}}\d u_1\ldots\d
  u_{k-1}}{ (u_1\ldots u_{k-1}(x-u_1-\ldots-u_{k-1}))^{3/2}}\, ,\qquad
x\geq 0\, ,
\end{equation}
the latter being $=0$ on $[0,k\eps]$. Let also
$$b(y)=\frac{1}{4\pi}\int_{y-\eps}^\eps\frac{\d u}{(u(y-u))^{3/2}}\,
,\qquad y\in(\eps,2\eps]\, .$$ Note that the functions $a_k,b$ depend
tacitly on the parameter $\eps$. Finally, we let
\begin{equation}\label{eq:7}
g(x)=\frac{\exp(-1/(4x))}{2\sqrt{\pi}x^{3/2}}\, ,\qquad x>0
\end{equation}
be the density of the hitting time of $1/\sqrt{2}$ by a standard
Brownian motion. This is the density of a stable law with parameter
$1/2$, which has Laplace transform
$$\int_{\R_+}\exp(-\lambda x)g(x)\d x=\exp(-\sqrt{\lambda})\, ,\qquad
\lambda\geq 0\, .$$ It can also be expressed as the density function
of the law of the total sum of the atoms of a Poisson random measure
on $\R_+$, with intensity
$$\frac{\d x}{2\sqrt{\pi}x^{3/2}}\ind_{\{x>0\}}\, .$$ The following
statement is a ``local limit theorem'' for $a$-skeletons of plane
trees.
\begin{prp}\label{sec:consequences-2}
  Let $\ms\in \T$.  Then for every $\eps\in(0,1)$ and $M>1/\eps$, the
  quantity
\begin{equation}\sup_{\substack{(\bx,\by)\in \Z_+^{\ms}\times \Z_+^{L(\ms)}:\\
    (\bx/\sqrt{n},\by/n) \in I_{M}(\ms)}}\!\!\!\! \left|
  n^{2|\ms|-1/2} \P_n({\rm Sk}_{\lf \eps n\rf}=
  (\ms,\bx,\by))-\psi^\ms(\bx/\sqrt{n},\by/n)\right|
\label{eq:loc-lim}
\end{equation}
converges to
$0$ as $n\to\infty$, where, setting $z=1-|\by|$,
$$\psi^\ms(\bx,\by)=\frac{\sqrt{\pi}}{2^{2|\ms|-2}}\left(
\sum_{k=0}^{\lf
  z/\eps\rf}\frac{(-|\bx|)^{k-2}}{k!}\int_{k\eps}^zg\left(
\frac{z-u}{|\bx|^2}\right)a_k(u)\d u\right)\prod_{v\in
  L(\ms)}b(y_v)\, .$$
\end{prp}

Although it is not apparent at first sight, the functions $\psi^\ms$
are equal to $0$ whenever $|\ms|>1/\eps$. Indeed, in this case $\ms$
cannot be the first component of the $\lf \eps n\rf$-skeleton of a
plane tree, since otherwise $\by_n$ would have strictly more than
$1/\eps$ components of size at least $\eps n$, which is impossible.

By comparing this proposition with a statement like Proposition
\ref{sec:brownian-crt-2}, it is tempting to interpret the functions
$\psi^\ms,\ms\in \T$ as density functions associated with a trimmed
version of the CRT. We are going to make this more precise in Section
\ref{sec:interpr-terms-psims}.

In order to prove the proposition, write
using Lemma \ref{sec:consequences-1},  
\begin{eqnarray*}
n^{2|\ms|-1/2} \P_n({\rm Sk}_{\lf \eps n\rf}=
  (\ms,\bx,\by))&=&  \frac{n^{-3/2}}{2^{2|\ms|-1}\mu_n}\\
&&\times 
n\,P\left(\sum_{i=1}^{|\bx|}X_i=n-|\by|,\, \max_{1\leq i\leq |\bx|}
X_i\leq \lf \eps n\rf\right)\\
&&\times \prod_{u\in L(\ms)}n^2P(X_1+X_2=y_u,X_1\vee X_2\leq \lf \eps n\rf).
\end{eqnarray*}
The asymptotics of the first term is $\sqrt{\pi}/2^{2|\ms|-2}$ by (\ref{eq:8}). 
The two following Lemmas give an uniform approximation of the two other terms.

\begin{lmm}\label{sec:consequences-4}
For any $M>1/\eps$, we have the convergence
$$\lim_{n\to\infty}\sup_{l/n\in
  [\eps+M^{-1},2\eps]}\left|n^2P(X_1+X_2=l,X_1\vee X_2\leq \lf \eps
n\rf)-b(l/n)\right|=0\, ,$$
\end{lmm}

\proof By definition, for $l/n\in [\eps+M^{-1},2\eps]$
\begin{eqnarray*}
\lefteqn{P(X_1+X_2=l,X_1\vee X_2\leq \lf \eps n\rf)=\sum_{j=l-\lf\eps
    n\rf}^{\lf \eps n\rf}\mu_j
  \mu_{l-j}}\\ &=&\frac{1}{4\pi}\sum_{j=l-\lf \eps n\rf}^{\lf \eps
  n\rf}\frac{1}{j^{3/2}(l-j)^{3/2}}+o\left(\sum_{j=\lf n/M\rf}^{\lf
  \eps n\rf}\frac{1}{j^{3/2}(l-j)^{3/2}}\right)\, ,
\end{eqnarray*}
as $n\to\infty$, where we have used the equivalent (\ref{eq:8}), and
the fact that $l-\lf \eps n\rf\geq \lf n/M\rf$. In turn,
\begin{eqnarray*}
\frac{1}{4\pi}\sum_{j=l-\lf \eps n\rf}^{\lf \eps
  n\rf}\frac{1}{j^{3/2}(l-j)^{3/2}}&=&\frac{1}{4\pi}\int_{l-\lf\eps
  n\rf}^{\lf \eps n\rf}\frac{\d x}{\lf x\rf^{3/2}(l-\lf x\rf)^{3/2}}\\
&=&\frac{1}{4\pi n^2}\int_{(l-\lf \eps n\rf)/n}^{\lf \eps
  n\rf/n}\frac{\d x}{(\lf nx \rf/n)^{3/2}(l/n-\lf nx\rf/n)^{3/2}}\, .
\end{eqnarray*}
This yields the result as the function $(x,y)\mapsto
x^{-3/2}(y-x)^{-3/2}$ is uniformly continuous on the compact set
$\{(x,y):(x,y-x)\in [(2M)^{-1},2\eps]^2\}$, in which the terms $(\lf
nx\rf/n,l/n)$ of the above integral are constrained to lie (at least
for $n$ large enough). \cq

\begin{lmm}\label{sec:consequences-5}
For $M>1$, 
\begin{eqnarray*}
\lefteqn{\lim_{n\to\infty}\sup_{\substack{l/\sqrt{n}\in[M^{-1},M]\\ m/n\in[0,1]}}
  \left|nP\left(\sum_{i=1}^{l}X_i=m,\max_{1\leq i\leq l}X_i\leq
  \lf\eps n\rf\right)\right.}\\ &&\hskip3cm-\left.\sum_{k\geq
  0}\frac{(-l/\sqrt{n})^{k-2}}{k!}\int_0^{m/n}\d u\, a_k(u)\,
g\left(\frac{m/n-u}{ (l/\sqrt{n})^2}\right)\right|=0\, .
\end{eqnarray*}
\end{lmm}

\proof For events $A,B_j,j\in J$ with finite $J$, one has
$$E\left[\ind_A\prod_{j\in J}(1-\ind_{B_j})\right]=\sum_{C\subseteq
  J}(-1)^{\#C}E\left[\ind_A\prod_{j\in C}\ind_{B_j}\right]\, ,$$ which
is usually called the {\em inclusion-exclusion principle}. Using this
and an elementary exchangeability argument, one has
\begin{eqnarray*}
n
P\left(\sum_{i=1}^{l}X_i=m,\max_{1\leq i\leq l}X_i\leq \lf\eps
n\rf\right)
=n\sum_{k=0}^{l}(-1)^k\binom{l}{k}P\left(\sum_{i=1}^{l}X_i=m,\min_{1\leq
  i\leq k}X_i>\lf \eps n\rf\right)\,;
\end{eqnarray*}
in order to let the reader follows more easily the step of the computations, we rewrite the RHS as
\begin{equation}\label{step-f}
\sum_{k=0}^{l}(-1)^k\frac{\binom{l}{k}}{n^{k/2}}\times n^{k/2+1} 
P\left(\sum_{i=1}^{l}X_i=m,\min_{1\leq   i\leq k}X_i>\lf \eps n\rf\right).
\end{equation}
Note that the sum can be reduced to its $m/\lf\eps n\rf$ first
terms, because the constraint on the minimum of the first $k$
variables $X_i$ forces $k\lf \eps n\rf<m$. Thus, there are at most
$1/\eps$ non-zero terms in the sum, which is independent on $n$. It
remains to show the (uniform) convergence of these individual terms. 

It is obvious that $\binom{l}{k}/n^{k/2}\sim (l/(\sqrt{n}))^k/k!$ uniformly for
$l/\sqrt{n}\in [M^{-1},M]$. To handle the probability term in (\ref{step-f}), we
introduce 
\begin{eqnarray*}
A_k(n,r)&=&P\left(X_1+\ldots+X_k=r,\min_{1\leq i\leq k}X_i>\lf\eps
n\rf\right)\\ G(l,r)&=&P(X_1+\ldots+X_l=r)\, ,
\end{eqnarray*}
and write the quantity after the $\times$ sign, in (\ref{step-f}),
\begin{eqnarray}\label{eq:10}
 \lefteqn{\sum_{r=0}^{m}n^{k/2+1}A_k(n,r)G(l-k,m-r) }\label{eq:15}\\
 &=&\int_0^{m/n}\!\left(n^{k/2+1}A_k(n,\lf rn\rf)\right)\left(\frac{\sqrt{n}}{l}\right)^2
 \left(l^2G(l-k,m-\lf rn\rf)\right)\d
  r\, . \nonumber
\label{eq:14}
\end{eqnarray}
In order to handle the terms involving $G(l,m)$, we use the
Gnedenko-Kolmogorov local limit theorem for sums of i.i.d.\ random
variables in the domain of attraction of a stable (1/2) distribution
\cite[Theorem 8.4.1]{BGT}, giving that
\begin{equation}\label{eq:16}
  \lim_{l\to\infty}\sup_{m\geq 0}\left|l^2G(l,m)-g(m/l^2)\right|=0\, .
\end{equation}
 In particular, in the expression (\ref{eq:10}), we can omit the
$-k$ inside the $G$ terms.  On the other hand, we have
\begin{equation}\label{eq:11}
\lim_{n\to\infty}\sup_{0\leq m\leq n}\left|n^{1+k/2}A_k(n,m)-
a_k(m/n)\right|=0\, .
\end{equation} Indeed, using (\ref{eq:8}),
\begin{eqnarray*}
A_k(n,m)&\sim& \sum_{\substack{j_1,\ldots,j_k>\lf\eps
    n\rf\\j_1+\ldots+j_k=m}}\prod_{l=1}^{k}\frac{1}{2\sqrt{\pi}j_l^{3/2}}\\ &=&
(4\pi)^{-k/2}\int \frac{\d v_1\ldots \d
  v_{k-1}\ind_{D(k,m)}(\bv)}{\left(m-\sum_{i=1}^{k-1}\lf
  v_i\rf\right)^{3/2}\prod_{i=1}^{k-1}\lf v_i\rf ^{3/2}}\, ,
\end{eqnarray*}
where by definition $D(k,m)=\{(v_1,\ldots,v_{k-1}):\min_{1\leq i\leq
  k-1}\lf v_i\rf\wedge (m-\sum_{i=1}^{k-1}\lf v_i\rf)>\lf \eps
n\rf\}$. After the linear change of variables $u_i=v_i/n$,
(\ref{eq:11}) boils down to a situation similar to Lemma
\ref{sec:consequences-4}, using the fact that the function
$$(u_1,\ldots,u_{k-1},x)\mapsto\left((x-u_1-\ldots-u_{k-1})
\prod_{i=1}^{k-1}u_i\right)^{-3/2}$$ is (uniformly) continuous on the
compact
$$\left\{(u_1,\ldots,u_{k-1},x):0\leq x\leq 1,\min_{1\leq i\leq k-1}
u_i\geq \eps/2,x-u_1\ldots-u_{k-1}\geq \eps/2\right\}\, .$$
It is now easy to conclude, by (\ref{eq:15}), (\ref{eq:16}) and
(\ref{eq:11}). \cq

\subsection{The case of trees}\label{sec:case-trees-3}

We now want to make a study parallel to the previous one in case of
trees. We let 
\begin{equation}\label{etet2}
\nu_n=\varrho^n\#\bT_n,n\geq 1,
\end{equation}
which defines a
probability distribution on $\N$ thanks to (\ref{eq:12}), satisfying 
\begin{equation}\label{eq:17}
\nu_n\build\sim_{n\to\infty}^{}\frac{\bc}{2\sqrt{\pi}n^{3/2}}\, .
\end{equation}
We now let $(\Omega,\FF,P)$ be a probability space under which an
i.i.d.\ sequence $X'_1,X'_2,\ldots$ of $\nu$-distributed random
variables is defined. A (non-plane) tree analog of Lemma
\ref{sec:consequences-1} states as follows. We let $G_a$ denote the
set of $a$-good trees, and $\bP_n^a:=\bP(\, \cdot\cap G_a)$.

\begin{lmm}\label{sec:case-trees-4}
  Let $(\ms,\bx,\by)$ be an $a$-good skeleton. Then, seeing
  ${\rm Sk}_a:\bt\mapsto {\rm Sk}_a(\bt)$ as a random variable,
\begin{eqnarray}
  \bP_n^a({\rm
    Sk}_a=(\ms,\bx,\by)) &=&\frac{1}{2^{|\ms|}\nu_n}
  P\left(\sum_{i=1}^{|\bx|}X'_i=n-|\by|\, ,\, \max_{1\leq i\leq |\bx|}
    X'_i\leq a\right)\nonumber\\
  & & \times\prod_{u\in L(\ms)}P(X'_1+X'_2=y_u,X'_1\vee
  X'_2\leq a,X'_1\neq X'_2).
  \nonumber
\end{eqnarray}
\end{lmm}

\proof 
The proof goes like that of Lemma \ref{sec:consequences-1}, applying
Lemma \ref{sec:case-trees-2} (and using $\nu$ instead of $\mu$); the only important difference
being that in the first displayed formula of the proof, we must
replace plane trees with trees, the factor $2^{|\bx|}$ does not appear
anymore, and we should consider pairs of trees $(\bt,\bt')$ such that
$|\bt|+|\bt'|=y_i$ and $|\bt|<|\bt'|$, the last constraint being
absent in the plane case. This can be replaced by the constraint
$|\bt|\neq |\bt'|$ by introducing factors of $2$, one for each $1\leq
j\leq |\ms|$.
\cq

\medskip

Besides the fact that the random variables $X'_i$ replace $X_i$, the
main difference with Lemma \ref{sec:consequences-1} is the presence of
the event $G_a$, and of the additional condition $X'_i\neq X'_2$ in
the last event. We can also make a comment on the fact that a
prefactor of the form $2^{-|\ms|}$ replaces a factor $2^{-2|\ms|+1}$
in Lemma \ref{sec:consequences-1}: if we were concerned with the law
of the non-plane version $p(\TT_n)$ of a $\P_n$-distributed random
variable $\TT_n$, then for each (non-plane) possible skeleton
$(\ms,\bx,\by)$ without symmetries, there would be exactly
$2^{|\ms|-1}$ plane skeletons contributing to its weight.

From this, one deduces a counterpart of Proposition
\ref{sec:consequences-2} for (non-plane) trees. For $\ms\in \T$, we
let $\tilde{I}_M(\ms)$ be defined in a similar way as $I_M(\ms)$ in
(\ref{eq:13}), except that we further ask that the components of $\bx$
are pairwise distinct, and are such that $x_{u1}>x_{u2}$ whenever
$u\in I(\ms)$.

\begin{prp}\label{sec:proof-theor-refs-1}
  Let $\ms\in \T$. Then for every $\eps\in(0,1)$ and $M>1/\eps$, the
  quantity
$$\sup_{\substack{\bx,\by\in 
    \Z_+^{\ms}\times \Z_+^{L(\ms)}:\\(\bx/\sqrt{n},\by/n) \in
    \tilde{I}_{M}(\ms)}}\!\!\!\!\! \left| n^{2|\ms|-1/2}
  \bP_n^{\lf \eps n\rf}({\rm Sk}_{\lf \eps n\rf}=
  (\ms,\bx,\by))-2^{|\ms|-1}\bc^{2|\ms|-1}\psi^\ms(\bc\bx/\sqrt{n},\by/n)\right|
$$
converges to $0$ as $n\to\infty$. 
\end{prp}

The proof follows exactly the same route as that of Proposition
\ref{sec:consequences-2}, using the two following intermediate lemmas.

\begin{lmm}\label{sec:case-trees-5}
For any $M>1/\eps$, we have the convergence
$$\lim_{n\to\infty}\sup_{l/n\in
  [\eps+M^{-1},2\eps]}\left|n^2P(X'_1+X'_2=l,X'_1\vee X'_2\leq \lf
  \eps n\rf,X'_1\neq X'_2)-\bc^2 b(l/n)\right|=0\, ,$$
\end{lmm}

\begin{lmm}\label{sec:case-trees-6}
For $M>1$, 
\begin{eqnarray*}
  \lefteqn{\lim_{n\to\infty}\sup_{\substack{l/\sqrt{n}\in[M^{-1},M]\\ m/n\in[0,1]}}
    \left|nP\left(\sum_{i=1}^{l}X'_i=m,\max_{1\leq i\leq l}X'_i\leq
        \lf\eps n\rf\right)\right.}\\ &&\hskip3cm-\left.\sum_{k\geq
      0}\frac{(-\bc l/\sqrt{n})^{k-2}}{k!}\int_0^{m/n}\d u\, a_k(u)\,
    g\left(\frac{m/n-u}{ (\bc l/\sqrt{n})^2}\right)\right|=0\, .
\end{eqnarray*}
\end{lmm}

In turn, the proofs of these lemmas are exactly the same as that of
Lemmas \ref{sec:consequences-4} and \ref{sec:consequences-5}. In the
proof of Lemma \ref{sec:case-trees-5}, one just has to be careful of
the extra $\bc^2$ term which comes in front of $b$, by using the
equivalent (\ref{eq:17}) instead of (\ref{eq:8}), and take into
account the (asymptotically unimportant) fact that $X'_1,X'_2$ are
constrained to take distinct values. In the proof of Lemma
\ref{sec:case-trees-6}, a similar remark holds for the function
$a_k$, which has to be replaced by $\bc^k a_k$, while the function $g$
must be replaced by $x\mapsto \bc^{-2}g(\bc^{-2}x)$, which is the
density probability function of the total sum of masses of a Poisson
process with intensity
$$\frac{\bc\,  \d x}{2\sqrt{\pi}x^{3/2}}\ind_{\{x>0\}}\, .$$

\section{Trimming continuous trees}\label{sec:density-psitt-its}

We now want to have an interpretation of the functions $\psi^\ms$ in
terms of the CRT. To this end, we first give another
interpretation of the trimming operations discussed in Section
\ref{sec:case-plane-trees}.

\subsection{Excursion functions and trimming}\label{sec:trimm-read-cont}

Recall that ${\cal E}$ is the set of continuous functions $f$ defined
on a compact interval $[0,\sigma_f]$, that are non-negative, and
satisfy $f(0)=f(\sigma_f)=0$.  With every $f\in {\cal E}$, one
associates an $\R$-tree $\TT_f=[0,\sigma_f]/\equiv_f$ endowed with the
distance $d_f$, as in Section \ref{sec:brown-cont-rand}. This tree is
naturally {\em rooted} (i.e.\ comes with a distinguished point) at the
point $\rho$, which is the $\equiv_f$-equivalence class of $0$.

For $I$ a real interval, we let $|I|=\sup I-\inf I$ denote its length.
For every (nonempty, open or closed) interval $I$, $a=\inf I$ and
$b=\sup I$, such that $\inf_I f= f(a)=f(b)$, we can define a new
excursion function $f_I\in {\cal E}$ by
$$f_I(t)=f(t+a)-f(a),0\leq t\leq
\sigma_{f_I}=|I|\, .$$ The extremities $a,b\in I$ encode the same
point $x=x_I$ in $\TT_f$, and $f_I$ encodes in turn an $\R$-tree
$\TT_f(x):=\TT_{f_I}$, which is naturally interpreted as the subtree
of $\TT_f$ rooted at $x$, i.e.\ the set of points $y\in \TT_f$ such
that
$$d_f(\rho,y)=d_f(\rho,x)+d_f(x,y)\, .$$
Yet otherwise said, this is the set of all descendants of $x$.  We
define the {\em mass} of $\TT_f(x)$ as $|\TT_f(x)|=|I|$. Finally, for
$a\in (0,\sigma_f)$ a positive real number, we let
$$\TT_f^{[a]}=\{x\in \TT_f:|\TT_f(x)|\geq a\}\, ,$$
which is also the closure of $\{x\in \TT_f:|\TT_f(x)|>a\}$, as is
easily checked. The set $\TT_f^{[a]}$ is an $\R$-tree which contains
$\rho$, and is naturally rooted at this point.

As we will soon see, the operation
$\TT_f\mapsto \TT_f^{[a]}$ is a continuous analog of the trimming
operation $\t\mapsto \t[a]$ defined in Section
\ref{sec:trimming-trees}.

First we present a slightly different way to look at $\TT_f^{[a]}$.
For any $f\in {\cal E}$ and $t\geq 0$, let $\O_f(t)$ be the set of
(non-empty) connected components of the open set $\{f>t\}$. It is
convenient to imagine $\O_f(t),t\geq 0$ as a {\em fragmentation
  process}, as the set $\{f>t\}$ decreases from $[0,\sigma_f]$ to
$\varnothing$ when $t$ moves from $0$ to $\infty$. Note that, if $I\in
\O_f(t)$, then $f_I$ is well-defined and is positive on
$(0,\sigma_{f_{I}})$, and as above, $I$ encodes a point $x_I$ in
$\TT_f$. This point will be in $\TT_f^{[a]}$ if and only if $|I|\geq
a$. This identifies $\TT_f^{[a]}$ with the set of intervals $I\in
{\cal O}_f(t)$ with $t\geq 0$ and $|I|\geq a$. Note however that this
identification is not one-to-one, as several intervals can encode the
same point of $\TT_f$.  

From this, we can state the following Lemma allowing to control the
Gromov-Hausdorff distance between $\TT_f$ and $\TT_f^{[a]}$.

\begin{lmm}\label{sec:excurs-funct-trimm-2}
  For every $f\in {\cal E}$ such that $\sigma_f>a$, it holds that
$$\dgh(\TT_f,\TT_f^{[a]})\leq \omega(f,a)\, ,$$  
where $\omega(f,a)=\sup_{|x-y|\leq a}|f(x)-f(y)|$ is the modulus of
continuity of $f$, evaluated at $a$. In particular, $\TT_f^{[a]}\to
\TT_f$ for the Gromov-Hausdorff metric as $a\to 0$.
\end{lmm}

\proof Since $\TT_f^{[a]}$ is defined as a subset of $\TT_f$, it
suffices to show that any point of $\TT_f$ not in $\TT_f^{[a]}$ is at
distance at most $\omega(f,a)$ from a point in $\TT_f$. So let $x$ be
such a point. For every $0<t<d_f(\rho,x)$, there is a unique interval
$I_x(t)\in {\cal O}_f(t)$ that contains some (hence all)
representative of $x$ in $[0,\sigma_f]$ for the equivalence relation
$\equiv_f$. Since $x$ is not in $\TT_f^{[a]}$ is holds that
$|I_x(t)|<a$ for $t$ close enough to $d_f(\rho,x)$, so let $t_0$ be
the infimum of all $t$'s with this property. The intervals
$I_x(t_0-\eps)$ decrease as $\eps\to 0$ to some interval $I_x(t_0-)$
with length $\geq a$, which encodes a vertex $x_0\in \TT_f^{[a]}$.
For every $\eps\in(0,d_f(\rho,x)-t_0)$, the interval $I_x(t_0+\eps)$
has length $<a$, so easily $d_f(x,x_{I_x(t_0+\eps)})=f(s)-f(s')\leq
\omega(f,a)$ where $s$ is a representative of $x$ and $s'$ is the left
end of $I_x(t_0+\eps)$. Hence $d_f(x,x_0)\leq \eps+\omega(f,a)$ for
every $\eps>0$, giving the result.  \cq

\subsection{The $a$-real skeleton description of
  $\TT_f^{[a]}$}\label{sec:encoding-tt_fa-as}

Let us now derive an alternative representation of $\TT_f^{[a]}$ in
terms of marked trees. 
To avoid trivialities we let
\[{\cal E}_a=\{f\in {\cal E}:\sigma_f>a\}\, ,\] and will assume from now on that $a>0$ and $f\in
{\cal E}_a$ so that in particular, $\TT_f^{[a]}$ is non-empty and not
reduced to a single point. 

\begin{figure}[htbp]
\psfrag{a}{$a$}
\psfrag{s}{$s_a(f)$}
\psfrag{t}{$t_a(f)$}
\psfrag{tt}{$t$}
\centerline{\includegraphics[height=2.3cm,width=14 cm]{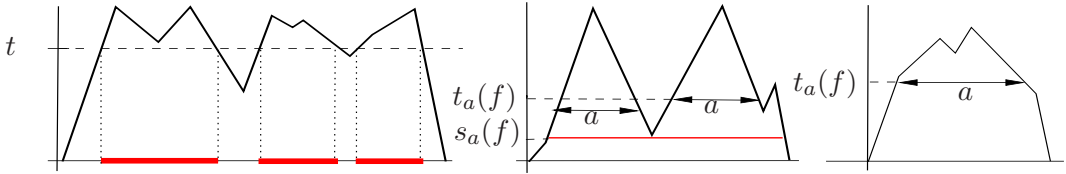}}
\caption{The set of intervals $\O_f(t)$ (left), the case
  $s_a(f)<+\infty$ (middle), the case $s_a(f)=+\infty$ (right).}
\label{frag}
\end{figure}

Assuming for a second that $\O_f(0)=(0,\sigma_f)$, by continuity of
$f$, $\O_f(t)$ contains a single interval of length $>a$ for every
small enough $t>0$. When $t$ increases, this property will break down
either because a second interval of length $>a$ is created, or if no
interval of length $>a$ remains (see Figure \ref{frag}). For every
$f\in{\cal E}_a$, we let
$$t_a(f)=\inf\left\{t\geq 0:\sup_{I\in \O_f(t)} |I|\leq a\right\} \in
 [0,\sup f]\, ,$$  and
$$s_a(f)=\inf\{t\geq 0:\exists I, I' \in \O_f(t):I\neq I',|I|>a
\mbox{ and }|I'|>a\}\, \in[0,\infty]\, .$$

By definition, for every $t<s_a(f)\wedge t_a(f)$, $\O_f(t)$ contains a
single interval $I_f(t)$ with $|I_f(t)|>a$, and these decrease to a
closed interval $I_f(s_a(f)\wedge t_a(f)-)$ as $t\uparrow s_a(f)\wedge
t_a(f)$. By convention, we let $I_f(0-)=(0,\sigma_f)$ in the case
$s_a(f)\wedge t_a(f)=0$.  It is immediate to see that if
$s_a(f)<\infty$ then $s_a(f)<t_a(f)$. In the case $s_a(f)=\infty$, we
let
$$Y_a(f):=|I_f(t_a(f)-)|\geq a\, .$$ In all cases, thanks to the
description of $\TT_f^{[a]}$ in terms of the intervals of ${\cal
  O}_f(t)$, one checks easily that the mapping $t\in[0,s_a(f)\wedge
  t_a(f))\mapsto x_{I_f(t)}\in\TT_f$ is an isometry, and that its
  image is exactly the set
$$\left\{x\in \TT_f^{[a]}:d_f(x,\rho)<s_a(f)\wedge t_a(f)\right\}\, .$$ This
isometry can be extended continuously to the whole segment
$[0,s_a(f)\wedge t_a(f)]$ by mapping its right end to
$x_{I_f(s_a(f)\wedge t_a(f)-)}$. Otherwise said, for $f\in {\cal
  E}_a$, the ``bottom'' of the tree $\TT_f^{[a]}$ is made of a line
segment (possibly reduced to the single point $\{\rho\}$). Moreover,
if $s_a(f)=\infty$, then the above isometry admits all of
$\TT_f^{[a]}$ as its image, so the latter is isometric to the segment
$[0,t_a(f)]$. So let us study what happens above level $s_a(f)$ in the
case $s_a(f)<\infty$.

We will be interested in a particular type of functions of ${\cal
  E}_a$.  We say that $x\in (0,\sigma_f)$ is a local minimum of $f$ if
$f(x)\leq f(y)$ for every $y$ in an open neighborhood of $x$. Note
that $0,\sigma_f$ are purposely excluded with this definition. We let
$m_f$ be the set of local minima of $f$.  

For $a>0$, we let $\ov{{\cal E}}_a$ be the set of functions $f\in
{\cal E}$ such that for all $x,y\in m_f, x\leq y$ such that
$f(x)=f(y)$, then $x=y$ or $\min\{ f(t)~:~ t \in [x,y]\}< f(x)$. This
condition implies that ${\cal T}_f$ has only binary branching points.
We also let $\ov{{\cal E}}_a^*$ be the subset of $\ov{{\cal E}}_a$
made of those functions $f$ such that
$$a\notin {\rm closure}(\{|I|:\exists t\geq 0,I\in \O_f(t)\})\, .$$
The sets $\ov{{\cal E}}_a,\ov{{\cal E}}_a^*$ are not nice from a
topological perspective (in particular they are neither open nor
closed), but it turns out that the mapping $f\mapsto \O_f$ behaves
well when restricted to these sets, as the following lemmas will show.

\begin{lmm}\label{sec:cont-funct-trimm-3}
  Let $a>0$ and $f\in \ov{{\cal E}}_a$.
  \begin{enumerate}
  \item If $s_a(f)=\infty$ and $f\in \ov{{\cal E}}_a^*$, then
    $Y_a(f)>a$. 
  \item If $s_a(f)<\infty$, then $\O_f(s_a(f))$ has exactly two
    components $I^{(l)}_f,I^{(r)}_f$ with $|I^{(l)}_f|\wedge
    |I^{(r)}_f|>a$, and these can be ordered so that $\sup
    I^{(l)}_f=\inf I^{(r)}_f$. We then set
    $f^{(l)}=f_{I^{(l)}_f},f^{(r)}=f_{I^{(r)}_f}$.
\end{enumerate}
\end{lmm}

\proof 1. Assume $s_a(f)=\infty$. The possibility that $Y_a(f)=a$ is
excluded by definition of $Y_a(f)$ and the fact that $a$ is not a
limit point of the set of lengths of intervals in $\cup_{t\geq
  0}\O_f(t)$. Hence $Y_a(f)>a$. 

2. Assume $s_a(f)<\infty$. Let $I_0=I_f(s_a(f)-)$. Necessarily, it
must hold that $f_{I_0}$ attains the value $0$ at some point $x$ of
$(0,\sigma_{f_{I_0}})$, which is unique as it must correspond to a
local minimum of $f\in \ov{{\cal E}}_a$. Hence, $\O_{f_{I_0}}(0)$ is
made of two intervals $I_1,I_2$ such that $\sup I_1=\inf I_2=x$, and
$\inf I_1=0,\sup I_2=|I_0|$. These two intervals correspond to two
intervals
$$I^{(l)}_f=(\inf I_0,x+\inf I_0)\, ,\qquad 
I^{(r)}_f=(x+\inf I_0,\sup I_0)\, ,$$ belonging to
$\O_f(s_a(f))$. By definition, the latter set contains two intervals
with lengths $>a$, but $I^{(l)}_f,I^{(r)}_f$ are the only possible
ones, since an interval $I\in \O_f(s_a(f))$ is a subinterval of some
interval in $\O_f(t)$ for $t<s_a(f)$, so in order that $|I|>a$ we must
have $I\subset I_f(t)$ for every $t<s_a(f)$, giving $I\subset
I_0$. \cq

\medskip

Point 2.\ in Lemma \ref{sec:cont-funct-trimm-3} entails that in the
case $s_a(f)<\infty$, the subtree of $\TT_f$ above the point
$x_{I_f(s_a(f)-)}$ splits into two branches with masses $>a$. More
precisely, if $h>0$ and $\TT,\TT'$ are two rooted $\R$-trees with
roots $\rho,\rho'$, we let $\Theta(h,\TT,\TT')$ be the $\R$-tree
obtained by identifying the roots $\rho,\rho'$ and grafting to this
point a line segment of length $h$. The latter can be easily
formalized as the quotient of the disjoint union
$$[0,h]\sqcup \TT\sqcup \TT'$$
in which the point $h\in [0,h]$ is identified with $\rho$ and $\rho'$,
and endowed with the proper distance function that extends the metrics
on each of these components.  The tree $\Theta(h,\TT,\TT')$ is
naturally rooted at the point $0\in [0,h]$. Then if $f\in \ov{{\cal
    E}}_a$ and $s_a(f)$, with the notations of Lemma
\ref{sec:cont-funct-trimm-3},
\begin{equation}
  \label{eq:20}
  \TT_f^{[a]}\build=_{}^{isom}
\Theta(s_a(f),\TT_{f^{(l)}}^{[a]},\TT_{f^{(r)}}^{[a]})\, ,
\end{equation}
while $\TT_f^{[a]}$ is isometric to $[0,t_a(f)]$ if $s_a(f)=\infty$.
Put together, these two facts allow a recursive construction of
$\TT_f^{[a]}$, which ends in a finite number of steps, since by
definition, we have $\sigma_{f^{(l)}}\vee \sigma_{f^{(r)}}\leq
\sigma_f -a$ whenever $s_a(f)<\infty$.

This description can be further simplified as follows.  Recall the
definition of the concatenation operation $\Xi$ for labeled plane
trees, defined in Section \ref{sec:case-plane-trees}.  We can
associate, in a recursive manner, a marked plane tree
$(\ms,\bg)=\zeta(f,a)$ with $\bg=(\bx,\by)\in \R^{\ms}\times
\R^{\ms}$ with every $f\in \ov{{\cal E}}_a$ in the following way.
\begin{itemize}
\item
If $s_a(f)=\infty$ then
$\ms=\{\varnothing\},x_\varnothing=t_f(a),y_\varnothing=Y_a(f)$. 
\item If $s_a(f)<\infty$ then
  $\zeta(f,a)=\Xi((s_a(f),0),\zeta(f^{(l)},a),\zeta(f^{(r)},a))$.
\end{itemize}
Note that every $\by$-mark of an internal vertex is $0$, so we simply
forget these marks, identifying $\zeta(f,a)$ with the tree $\ms$ with
marks $(\bx,\by)\in \R^{\ms}\times \R^{L(\ms)}$.

\begin{figure}[htbp]
\psfrag{sv}{$s_\varnothing$}
\psfrag{s1}{$s_1$}
\psfrag{s2}{$s_2$}
\psfrag{t0}{$t_0$}
\psfrag{t1}{$t_1$}
\psfrag{t2}{$t_2$}
\psfrag{Y1}{$Y_1$}
\psfrag{Y2}{$Y_2$}
\psfrag{Y3}{$Y_3$}

\psfrag{}{}
\psfrag{}{}
\psfrag{}{}
\psfrag{}{}
\psfrag{}{}
\centerline{\includegraphics[height=4cm]{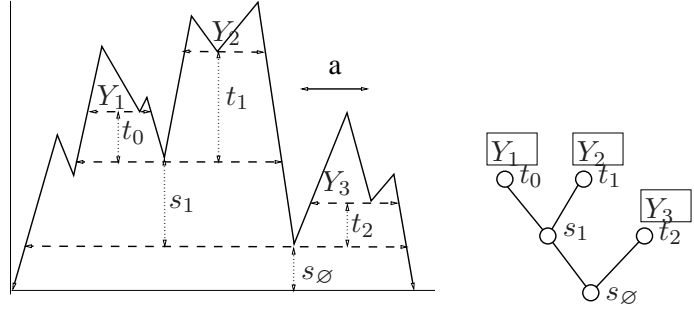}}
\caption{Representation of $\zeta(f,a)$}
\label{dec}
\end{figure}

The marked tree $\zeta(f,a)$ is called the $a$-real skeleton of $f$
(or of $\TT_f$). Building the tree $\TT_f^{[a]}$ from $\zeta(f,a)$ is
now an easy matter by comparing the definition of $\zeta(f,a)$ with
(\ref{eq:20}). With every pair of the form $(\ms,\bx)$, where
$\ms\in\T$ and $\bx=(x_u,u\in \ms)\in (0,\infty)^\ms$, we associate an
$\R$-tree $\theta(\ms,\bx)$ in a natural way, by grafting segments
with lengths $x_u$ according to the genealogy of $\ms$. Formally, 
we let $\theta((\varnothing,x_\varnothing))$ be the isometry class of
$[0,x_\varnothing]$, and recursively, if $\ms\neq \{\varnothing\}$, 
$$\theta(\ms,\bx)=\Theta(x_\varnothing,\theta((\ms,\bx)_{1}),
\theta((\ms,\bx)_{2}))\, ,$$ where by definition $(\ms,\bx)_{u}$ is
the plane subtree $\ms_u$ with labels $\bx_u=(x_{uv},v\in \ms_u)$.

\begin{lmm}\label{sec:a-real-skeleton}
  Let $f\in \ov{{\cal E}}_a$.  Then $\TT_f^{[a]}$ and
  $\theta(\zeta(f,a))$ are isometric, where in the latter notation we
  do not take the $\by$-marks into account.
\end{lmm}

To prove this, observe that $\theta(\zeta(f,a))$ is isometric to
$[0,t_a(f)]$ if $s_a=\infty$, and that
$$\theta(\zeta(f,a))=
\Theta(s_a(f),\theta(\zeta(f^{(l)},a)),\theta(\zeta(f^{(r)},a)))\, ,$$
if $s_a(f)<\infty$, and compare with (\ref{eq:20}).

\subsection{Continuity properties of the
  $a$-real skeleton}\label{sec:cont-prop-a}

\begin{lmm}\label{sec:cont-funct-trimm}
  Let $a>0$ and $(a_n,n\geq 1)$ be a positive sequence converging to
  $a$. Let $(f_n,n\geq 1)$ be a sequence such that $f_n\in {\cal
    E}_{a_n}$ for every $n$, uniformly converging to $f\in \ov{{\cal E}}_a^*$. Then
\begin{enumerate}
\item
If $s_a(f)=\infty$, then $t_{a_n}(f_n)\to t_a(f)$, and
$Y_{a_n}(f_n)\to Y_a(f)$ as $n\to\infty$.
\item
If $s_a(f)<\infty$, then $s_{a_n}(f_n)\to s_a(f)$ and
$(f_n^{(l)},f_n^{(r)})\to (f^{(l)},f^{(r)})$ as $n\to\infty$.
\end{enumerate}
\end{lmm}

\proof 1. This comes immediately from the uniform convergence of
$f_n\to f$, and does not use the fact that $f\in \ov{{\cal E}}^*_a$.

2. Let $x=\inf I^{(l)}_f,y=\sup I^{(l)}_f=\inf I^{(r)}_f$ and $z=\sup
I^{(r)}_f$. Let also $\eps\in(0, ((y-x)\wedge (z-y))/2)$. From the
uniform convergence of $f_n$ to $f$, the minimum of $f_n$ over
$(x+\eps,z-\eps)$ is attained at some $y_n$, which converges to $y$ as
$n\to\infty$ since $f(v)>f(y)$ for all $v\in [x+\eps,y)\cup
(y,z-\eps]$. Let $x_n=\sup\{v<y_n:f_n(v)=f_n(y_n)\}$ and
$z_n=\inf\{v>y_n: f_n(v)=f_n(y_n)\}$. Then we have $x_n\to x,z_n\to z$
as $n\to\infty$. Consequently, the intervals $(x_n,y_n),(y_n,z_n)$
have both lengths $>a_n$ for $n$ large enough, since $y-x,z-y$ are
$>a$ and $a_n\to a$. It remains to note that for large $n$, there can
be no two intervals of length $>a_n$ in $\O_{f_n}(t)$ for
$t<f_n(y_n)$, in which case it will follow that
$s_{a_n}(f_n)=f_n(t_n)\to f(t)=s_a(f)$, and the remaining properties
are easy. But assuming that two such intervals exist, one of them must
contain the interval $(x_n,z_n)$. By extracting subsequences, we may
assume that the other interval is of the form $(x'_n,z'_n)$ (with
$z'_n-x'_n>a_n$) and that its extremities converge to $(x',z')$. By
construction $f(x')=f(z')\leq s_a(f)$ and by definition, it holds that
$a=\lim a_n\leq z'-x'\leq a$. Hence, $(x',z')$ is an interval of
length $a$ in $\O_f(f(x'))$ contradicting the fact that $f\in
\ov{{\cal E}}^*_a$.  \cq

\begin{crl}\label{sec:cont-funct-trimm-1}
  Let $a>0$ and $(a_n,n\geq 1)$ be a positive sequence converging to
  $a$. Let $(f_n,n\geq 1)$ be a sequence such that $f_n\in \ov{{\cal
      E}}_{a_n}$ for every $n$, converging to $f\in \ov{{\cal
      E}}_a^*$. Then $\zeta(f_n,a_n)\to \zeta(f,a)$. Consequently, it
  holds that $\dgh(\TT_{f_n}^{[a_n]},\TT_f^{[a]})\to 0$. 
\end{crl} 

The first part of this statement is obtained directly by combining
Lemma \ref{sec:cont-funct-trimm} with the recursive definition of
$\zeta(f,a)$. The last sentence follows from Lemma
\ref{sec:a-real-skeleton} and the easy fact that $\bx\mapsto
\theta(\ms,\bx)$ is continuous from $(0,\infty)^\ms$ to
$(\mathscr{T},\dgh)$, for every $\ms\in \T$.

\subsection{The case of discrete
  trees}\label{sec:case-discrete-trees}

\begin{figure}[htbp]
\centerline{\includegraphics[height=2cm]{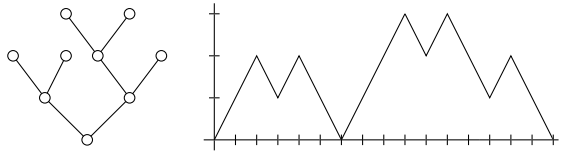}}
\caption{A plane tree and its contour process.}
\label{PF}
\end{figure}

We now reconnect the previous discussion to our study of (discrete)
trees. Let $\t\in \T_n$. The {\em contour exploration} of $\t$ is the
sequence $(u_\t(0),\ldots, u_\t(4n-4))$ where $u_\t(0)=\varnothing$,
and recursively, $u_\t(i+1)$ is the first child (in lexicographical
order) of $u_\t(i)$ not in $\{u_\t(0),\ldots,u_\t(i)\}$, if such a
child exists, or the parent of $u_\t(i)$ otherwise. Each vertex is
visited as many times as its degree in this sequence, which explains
that the procedure stops at the $(4n-4)$-th step (indeed, each edge of
the tree is visited twice, and there are $2n-2$ edges in a binary tree
with $n$ leaves). A picture of this is to imagine a particle wrapping
around the tree and moving at speed $1$, starting from $\varnothing$
(see Figure \ref{PF}).

\begin{defn}\label{sec:trimm-read-cont-1}
The contour process of $\t\in \T_n$ is the continuous function
$C_\t(t),0\leq t\leq 4n-4$ whose graph is obtained by a linear
interpolation between the points $(i,|u_\t(i)|),0\leq i\leq 4n-4$. 
\end{defn}

The contour process of $\t$ is also called the {\em Dyck path
  encoding} of $\t$, as the latter can be recovered from its contour
process. We will only be interested in recovering the metric space
$(\t,d_\t)$ from $C_\t$. If $v,w\in \t$, let $0\leq i,j\leq 4n-4$ be
such that $u_\t(i)=v,u_\t(j)=w$. Then it is easy to see that $|v\wedge
w|=\min_{i\wedge j\leq k\leq i\vee j}C_\t(k)$, independently of the
choice of $i,j$. Therefore,
\[d_\t(v,w)=C_\t(i)+C_\t(j)-2\min_{i\wedge j\leq k\leq i\vee
  j}C_\t(k)=d_{C_\t}(i,j)\,.\] Of course, for integers $i$ and $j$,
$d_{C_\t}(i,j)=0$ if and only if $u_\t(i)=u_\t(j)$, which is
equivalent to $i\equiv_{C_\t}j$, and so, the quotient space
$\{0,1,\ldots,4n-4\}/\equiv_{C_\t}$ endowed with the distance
$d_{C_\t}$, is isometric to $(\t,d_\t)$. Therefore $\t$ is isometric
to a subset of $\TT_{C_\t}$, and it is immediate to check that
$$  \dgh(\t,\TT_{C_\t})\leq 1\, .$$

\begin{figure}[htbp]
\psfrag{0}{0}
\psfrag{1}{1}
\psfrag{2}{2}\psfrag{3}{3}
\psfrag{4}{4}
\psfrag{5}{5}
\psfrag{12}{12}
\psfrag{16}{16}

\centerline{\includegraphics[height=7cm]{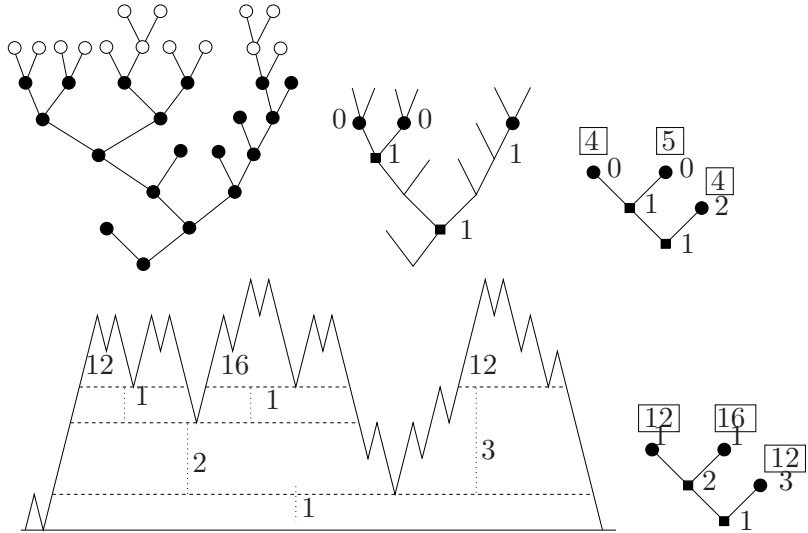}}
\caption{Comparison between the encoding via the skeleton and via the
  function $\zeta$. On the first line the trimmed tree $\t[3]$ and the
  skeleton representation. On the second line, representation of
  $\zeta(C_\t,12)$.}
\label{PF2}
\end{figure}

The trimming operation discussed in Section \ref{sec:case-plane-trees}
has a direct interpretation in terms of trees of the form
$\zeta(f,a)$, when $f$ is the contour function of a plane tree.

\begin{lmm}\label{sec:excurs-funct-trimm-1}
  Let $a>0$ be an integer, and $\t\in \T_n$. If
${\rm Sk}_a(\t)=(\ms,\bx,\by) $ then 
$$\zeta\left(C_{\t},4(a-1)\right)\,=(\ms,\bx',\by')\, ,$$ where
  $x'_u=x_u+ \ind_{\{u\neq \varnothing\}}$, for every $u\in \ms$, and
  $y'_u=4y_u-4$ for any $u\in L(\ms)$.
\end{lmm}
 Figure \ref{PF2} illustrates precisely the Lemma and the differences
 between the two points of view: the labels $\bx$ of the vertices are
 the same up to 1 due to discrete-continuous artifacts, and a subtree
 of size $k$ becomes a subpath of length $4(k-1)$; the proof below
 formalises these observations.\medskip

\proof We first note that since the plane trees that we consider are
binary, it is immediate that the contour function $C_\t$ is an element
of $\ov{{\cal E}}_a$ for $a<4|\t|-4$, so that $\zeta(C_\t,a)$ is
well-defined. 

We will use the following alternative description of ${\rm
  Sk}_a(\t)=(\ms,\bx,\by)$.  Recall the definition of
$\pi_{\t[a]}:\ms\to I_2(\t[a])$. Then for $u\in \ms$ and $i\in\{1,2\}$
such that $ui\in \ms$, we have
$$x_{ui}=|\pi_{\t[a]}(ui)|-|\pi_{\t[a]}(u)|-1\, ,$$
counting the number of elements of $I_1(\t[a])$ comprised between
$\pi_{\t[a]}(u)$ and $\pi_{\t[a]}(ui)$, while
$x_\varnothing=|\pi_{\t[a]}(\varnothing)|$.

Let us compare this with the construction of $\zeta(C_\t,4(a-1))$. For
$t\geq 0$ an integer, an interval $I$ of ${\cal O}_{C_\t}(t)$ encodes
a vertex $u_I$ of $\t$, seen as a subset of $\TT_{C_\t}$ as in the
construction above. If the length of the interval is $\geq 4a-4$, then
the subtree $\t_{u_I}$ has at least $a$ leaves, and consequently, this
vertex belongs to $\t[a]$. On the other hand, if the length is $\leq
4a-4$, then none of the strict descendants of $u_I$ can be in $\t[a]$.
From these observations, it follows that in the case
$s_{4(a-1)}(C_\t)=\infty$, it holds that
$${\rm
  Sk}_a(\t)=(\{\varnothing\},s_{4(a-1)}(C_\t),Y_{4(a-1)}(C_\t)/(4(a-1)))\,
.$$ In the case $s_{4(a-1)}(C_\t)<\infty$, the vertex of $\t$ encoded
by the interval $I_{C_\t}(s_{4(a-1)}(C_\t)-)$ is the vertex
$\pi_{\t[a]}(\varnothing)$ of Section \ref{tro}, and it holds that
${\rm Sk}_a(\t)=(\ms,\bx,\by)$ satisfies $x_\varnothing=s_a(C_\t)$.
Now, it must be noted that the functions $C_\t^{(l)},C_\t^{(r)}$
defined as in Lemma \ref{sec:cont-funct-trimm-3} are not the contour
functions of a binary plane tree, because they are strictly positive
on their respective intervals of definitions. We let
$\t^{(l)},\t^{(r)}$ be the binary plane trees that admit as contour
functions the paths $C_\t^{(l)},C_\t^{(r)}$, to which the first and
last steps have been removed. 
Then
$${\rm Sk}_a(\t)=\Xi(s_{4(a-1)}(C_\t),{\rm Sk}_a(\t^{(l)}),
{\rm Sk}_a(\t^{(r)}))\, .$$ Using this and an inductive argument, this
indeed shows that
$$x_u=x'_u-1\, ,\qquad u\neq \varnothing\, ,$$ with the notations of
the statement, this $-1$ being due to the fact that the quantities
$s_a,t_a$ associated with the functions $C_\t^{(l)},C_\t^{(r)}$, must
be subtracted $1$ in order to yield the number of skeleton leaves
comprised between two elements of $I_2(\t[a])$. \cq

\medskip

In view of deriving scaling limits, define the rescaled contour
function of $\t$ by
$$c_\t(s)=\frac{C_\t(4(|\t|-1)s)}{\sqrt{2|\t|}}\, ,\qquad 0\leq s\leq
1\, ,$$ and the $n$-rescaled skeleton as
$${\rm Sk}^{(n)}_a(\t)=(\ms,n^{-1/2}\bx,n^{-1}\by)\, ,$$ where
$(\ms,\bx,\by)={\rm Sk}_a(\t)$. Then by Lemma
\ref{sec:excurs-funct-trimm-1}, it holds that for $\t\in \T_n$,
$$\zeta\left(c_{\t},\frac{4(a-1)}{4(n-1)}\right)\,=(\ms,\bx',n^{-1}\by)
,$$ where for any $u\in \ms$, $|n^{-1/2}\bx_u-\bx'_u|\leq n^{-1/2}$.

\subsection{Scaling limits}\label{sec:scaling-limits}

We now discuss the scaling limits of the trimmed tree. The core result
is the counterpart of Proposition \ref{sec:brownian-crt-2} (in fact,
this is a stronger result) in terms of contour processes; it is due to
Aldous \cite{aldouscrt93} (see also \cite{mm01})
\begin{prp}\label{sec:trimm-read-cont-2}
  Let $\TT_n,n\geq 1$ be random variables with respective
  distributions $\P_n$. Then $c_{\TT_n}$ converges to $2\ee$, $\ee$
  the normalized Brownian excursion, in distribution for the topology
  of uniform convergence on ${\cal C}([0,1])$.
\end{prp}

\begin{crl}\label{sec:cont-funct-trimm-2}
  For every $\eps\in(0,1)$, it holds that ${\rm Sk}^{(n)}_{\lf \eps
    n\rf}(\t)$ under $\P_n$ converges in distribution to
  $\zeta(2^{3/2}\ee,\eps)$. 
\end{crl}

Note the normalization in front of the Brownian excursion in the last
statement, which comes from the fact that the $\bx$-components of the
rescaled skeleton are obtained by dividing by $\sqrt{n}$ the
components of the skeleton (instead of $\sqrt{2n}$).  From Corollary
\ref{sec:cont-funct-trimm-1} and the discussion after Lemma
\ref{sec:excurs-funct-trimm-1}, this is an immediate consequence of
the following

\begin{lmm}\label{sec:excurs-funct-trimm}
For every $\eps\in(0,1)$, it holds that with probability $1$, $\ee$ is
an element of $\ov{{\cal E}}_\eps^*$. 
\end{lmm}

\proof The fact that the local minima are realized once (a.s.) is a
well-known property of the Brownian excursion, which is inherited from
the analogous property for Brownian motion, and the property of
positivity of $\ee$ over $(0,1)$ comes from the definition. Hence a.s. $\ee$
belongs to $\ov{{\cal E}}_\eps$. It remains to show that $\eps$ is not a limit point of the
lengths of intervals of $\cup_{t\geq 0}\O_\ee(t)$. To show this, we
use the fact from \cite{bertsfrag02} that if $U$ is a uniform random
variable in $[0,1]$, independent of $\ee$, then for $0\leq t<\ee_U$,
the unique element $I_U(t)$ of $\O_\ee(t)$ containing $U$ is such that
$(|I_U(t)|,t\geq 0)$ has same distribution as
$(\exp(-\xi_{\tau(t)}),t\geq 0)$, where
\begin{itemize}
\item
$\xi$ is a subordinator \cite{bertlev96}, with L\'evy measure $e^x\d
  x/(e^x-1)^{3/2}$, and
\item
$\tau$ is the inverse of the increasing, continuous process $t\mapsto
  \int_0^t\exp(-\xi_s/2)\d s$.
\end{itemize}
In particular, the closure $C$ of the range of $(|I_U(t)|,t\geq 0)$ is
the image of the closure of the range of $\xi$ by $x\mapsto e^{-x}$,
to which has been adjoined the point $0$. By well-known properties of
subordinators \cite{bertlev96}, points are polar for $\xi$, so that
$\eps$ is not in $C$, almost-surely. Since $U$ is uniform and
independent of $\ee$, the wanted property easily follows --- for
instance, one can note that a.s.\ $\cup_{t\geq 0}\O_\ee(t)$ is the set
of intervals $I_{U_i}(t),i\geq 1,0\leq t<\ee_{U_i}$ where $U_i,i\geq 1$
is i.i.d.\ uniform on $[0,1]$ and independent of $\ee$, and with the
obvious notation for $I_{U_i}(t)$. \cq

\subsection{Interpretation in terms of 
$\psi^\ms$}\label{sec:interpr-terms-psims}

Now we reconnect the previous discussion to the functions $\psi^\ms$
defined in Section \ref{sec:consequences}.

\begin{lmm}\label{sec:consequences-3}
  For every $\eps\in (0,1)$, the functions $\psi^\ms,\ms\in \T$ form a
  ``probability density'', in the sense that
$$\sum_{\ms\in \T}\int_{\R_+^\ms}\d \bx\int_{\R_+^{L(\ms)}}\d \by\,
\psi^\ms(\bx,\by)=1\, .$$ 
For every $\eps\in(0,1)$, we can define a probability
measure $\Psi(\d(\ms,\bx,\by))$ on $\{(\ms,\bx,\by):\ms\in \T,\bx\in
\R^{\ms},\by\in \R^{L(\ms)}\}$ by the formula
$$\langle \Psi,F\rangle=\sum_{\ms\in \T}\int_{\R_+^\ms}\d
\bx\int_{\R_+^{L(\ms)}}\d \by\, \psi^\ms(\bx,\by)F(\ms,\bx,\by)\, ,$$
for every non-negative measurable function $F$.  Then $\Psi$ is the
law of $\zeta(2^{3/2}\ee,\eps)$.
\end{lmm}

\proof Fix $\ms\in \T$. Recall the definition (\ref{eq:13}) of the
open set $I_M(\ms)\subset \R_+^\ms\times \R_+^{L(\ms)}$. Then,
defining $\lf\bx\rf$ by taking integer parts componentwise, we have
\begin{eqnarray*}
  \lefteqn{ \P_n\left({\rm Sk}^{(n)}_{\lf \eps n\rf}\in
    \sint\{(\ms,\bx,\by): (\bx,\by)\in I_M(\ms)\}\right)} \\&=&
  \sum_{\ms\in \T}\int_{I_M(\ms)}\d\bx \d\by \, n^{2|\ms|-1/2}\,
  \P_n({\rm Sk}_{\lf \eps n\rf}=(\ms,\lf \bx\sqrt{n}\rf,\lf
  n\by\rf))\, .
\end{eqnarray*}
The $\liminf$ of the left-hand side as $n\to\infty$ is at least
$P(\zeta(2^{3/2}\ee,\eps)\in \sint\{(\ms,\bx,\by): (\bx,\by)\in
I_M(\ms))$, because of Corollary \ref{sec:cont-funct-trimm-2} and
well-known properties of weak convergence of probability measures. And
the right-hand side converges to
$$ \sum_{\ms\in \T}\int_{I_M(\ms)}\d\bx \d\by \, \psi^\ms(\bx,\by)\,
,$$ because of the uniform convergence stated in Proposition
\ref{sec:consequences-2}.  This yields
$$P\left(\zeta(2^{3/2}\ee,\eps)\in \sint\{(\ms,\bx,\by): (\bx,\by)\in
I_M(\ms)\}\right) \leq \sum_{\ms\in \T}\int_{I_M(\ms)}\d\bx \d\by \,
\psi^\ms(\bx,\by)\, . $$ Letting $M\to \infty$, the sets $I_M(\ms)$
respectively increase to $\{(\bx,\by)\in
(0,\infty)^\ms\times(\eps,2\eps)^{L(\ms)}\}$, and we know from Lemma
\ref{sec:excurs-funct-trimm} that the law of $\zeta(2^{3/2}\ee,\eps)$
is supported on the union of such sets as $\ms$ ranges in
$\T$. Therefore, the probability on the left-hand side converges to
$1$ as $M\to\infty$, which yields the first statement.  The fact that
$\Psi$ is the law of $\zeta(2^{3/2}\ee,\eps)$ is a simple adaptation
of the previous argument, and is omitted. \cq

\section{Proof of Theorem \ref{sec:bc1}}\label{sec:proof-theor-refs-3}

We now finally embark in the proof of our main results. This will be
done in two steps: first, we prove that the $\lf\eps n\rf$-trimmed
tree obtained from a $\bP_n$-distributed tree converges to the
$\eps$-trimming of the (appropriately scaled) CRT. Then we show that
the trimmed tree is ``not too far'' from the original tree.

\subsection{Convergence of the $a$-skeleton under
  $\bP_n$}\label{sec:conv-a-skel}

The main result of this section is

\begin{prp}\label{sec:interpr-terms-psims-1}
  For every $\eps\in(0,1)$, the random variables
$$\bt\mapsto \frac{\bc}{\sqrt{2n}}\bt[\lf \eps n \rf]\, ,\qquad
\mbox{under }\bP_n$$ converge to $\theta(\zeta(2\ee,\eps))$ as $n\to\infty$,
in distribution for the topology of $(\mathscr{M},\dgh)$, where $\ee$
is the standard Brownian excursion.
\end{prp}

To prove this, we study the scaling limit of $\bt\mapsto {\rm Sk}_{\lf
  \eps n\rf}(\bt)$ under $\bP_n$. First, we define a probability
measure, by
$$\langle \Psi^\circ,F\rangle= \sum_{\ms\in
  \T}2^{|\ms|-1}\int_{\Delta_\ms}\d \bx\int_{\R_+^{L(\ms)}}\d \by\,
\bc^{2|\ms|-1}\psi^\ms(\bc\bx,\by)F(\ms,\bx,\by)\, ,$$ where
$$\Delta_\ms=\{(x_u,u\in \ms):x_{u1}>x_{u2} \mbox{ for all }u\in
I(\ms)\}\, .$$ The fact that $\langle \Psi^\circ,1\rangle=1$ is a
consequence of the fact that $\Psi^\circ$ charges only good labeled
trees $(\ms,\bx,\by)$, and Lemma \ref{sec:case-trees-1}, because of
the presence of the factor $2^{|\ms|-1}$. The following is the analog
of Corollary \ref{sec:cont-funct-trimm-2} (and the end of Lemma
\ref{sec:consequences-3}) in the context of (non-plane) trees.

\begin{lmm}\label{sec:interpr-terms-psims-3}
  For every $\eps\in(0,1)$, the law of $\bt\mapsto {\rm Sk}^{(n)}_{\lf
    \eps n\rf}(\bt)$ under $\bP_n$ converges weakly to the measure
  $\Psi^\circ$. 
\end{lmm}

\proof 
For $\ms\in \T$ and $\bx\in (0,\infty)^\ms,\by\in (0,\infty)^{L(\ms)}$, we set 
$$g_n(\ms,\bx,\by)=n^{2|\ms|-1/2}\bP_n^{\lf \eps n\rf}({\rm Sk}_{\lf
  \eps n\rf}=(\ms,\lf \bx\sqrt{n}\rf,\lf n\by\rf))\, .$$ Using
Proposition \ref{sec:proof-theor-refs-1}, we have that
$g_n(\ms,\bx,\by)\to 2^{|\ms|-1}\bc^{2|\ms|-1}\psi^\ms(\bc\bx,\by)$
pointwise, and uniformly over $\bigcup_\ms\tilde{I}_M(\ms)$ for every
$M>0$.  Consequently, for every $M>0$, it holds that
\begin{eqnarray*}
 \bP_n(G_{\lf \eps n\rf})&=&\sum_{\ms\in \T}\int_{\R_+^\ms\times
   \R_+^{L(\ms)}}\d\bx \d\by g_n(\ms,\bx,\by)\\&\geq& \sum_{\ms\in
   \T}\int_{\tilde{I}_M(\ms)}\d\bx \d\by g_n(\ms,\bx,\by)\\ &\build
 \longrightarrow_{n\to\infty}^{} & \sum_{\ms\in
   \T}2^{|\ms|-1}\int_{\tilde{I}_M(\ms)}\d\bx \d \by
 \bc^{2|\ms|-1}\psi^\ms(\bc\bx,\by)\, ,
\end{eqnarray*}
which as $M\to\infty$ converges to 
$$\sum_{\ms\in
  \T}2^{|\ms|-1}\int_{\Delta_\ms}\d\bx\int_{(\eps,2\eps)^{L(\ms)}} \d
\by \bc^{2|\ms|-1}\psi^\ms(\bc\bx,\by)=\langle \Psi^\circ,1\rangle=1\,
.$$ Therefore, $\lim\bP_n(G_{\lf \eps n\rf})=1$, and by the so-called
Scheffé Lemma, it holds that $g_n$ converges in $L^1$. Letting $F$ be a
uniformly continuous function bounded by $K$, it holds that
$$|\bP_n(F({\rm Sk}_{\lf \eps n\rf}^{(n)}))-\bP_n^{\lf \eps
  n\rf}(F({\rm Sk}_{\lf \eps n\rf}^{(n)}))|\leq K(1-\bP_n(G_{\lf \eps
  n\rf}))\to 0\, ,$$ while
$$\bP_n^{\lf \eps n\rf}(F({\rm Sk}_{\lf \eps
  n\rf}^{(n)}))=\sum_{\ms\in \T}\int_{\R_+^\ms\times
  \R_+^{L(\ms)}}\d\bx\d\by F(\ms,\lf \bx\sqrt{n}\rf/\sqrt{n},\lf
n\by\rf/n )g_n(\ms,\bx,\by)$$ converges to $\langle
\Psi^\circ,F\rangle$, using the uniform continuity of $F$ and the
$L^1$-convergence of $g_n$.  \cq

\medskip

\begin{lmm}
  If $(\ms,\bx,\by)$ and $(\ms',\bx',\by')$ are respectively $\Psi$
  and $\Psi^\circ$ distributed, then $\theta(\ms,\bx)/\bc$ and
  $\theta(\ms',\bx')$ have the same distribution.
\end{lmm}

\proof Note that $\theta(\ms,\bx)=\theta(\ms',\bx')$ whenever
$(\ms,\bx)\approx(\ms',\bx')$, i.e.\ whenever these two plane labeled
trees represent the same labeled tree. Therefore, for a given
$a$-good tree $\bt$, $\theta({\rm Sk}_a(\bt))=\theta(\ms',\bx')$ for
any of the $2^{|\ms|-1}$ choices of $(\ms',\bx')\approx {\rm
  Sk}_a(\bt)$, and there are exactly $2^{|\ms|-1}$ of them. The
conclusion follows easily from the definition of $\Psi^\circ$ and the
fact that $\theta$ is a continuous function.  \cq

\medskip

\noindent{\bf Proof of Proposition \ref{sec:interpr-terms-psims-1}. }
As a consequence of Lemma \ref{sec:interpr-terms-psims-3}, and the
easily checked fact that
$$\dgh(\theta({\rm Sk}_{\lf \eps n\rf}(\bt)),\bt[\lf \eps
  n\rf])\leq 1\, $$ for every $\bt\in \bT_n$, we obtain that
$n^{-1/2}\bt[\lf \eps n\rf]$ converges in distribution to the
image by $\theta$ of a random variable with distribution $\Psi^\circ$,
which by Lemma \ref{sec:interpr-terms-psims-3} and Lemma
\ref{sec:interpr-terms-psims-1} has same law as
$\theta(\zeta(2^{3/2}\ee/\bc,\eps))$. \cq

\subsection{Tightness}
\label{sec:proof-theor-refs}

By Proposition \ref{sec:interpr-terms-psims-1}, we know that the
$\lf\eps n\rf$-trimming of a $\bP_n$-distributed tree is close in
distribution, when $n$ is large, to $\theta(\zeta(2\ee,\eps))$. In
turn, the latter is close to the CRT $\TT=\TT_{2\ee}$ when $\eps$ is
small, by Lemma \ref{sec:excurs-funct-trimm-2}. On an intuitive basis,
the proof will be complete if we are able to show that for $\eps$
small enough, under $\bP_n$, trees $\bt$ are typically close, in the
Gromov-Hausdorff sense, to $\bt[\lf \eps n\rf]$ in a uniform way as
$n\to\infty$.  This is what the next lemma is taking care of.

\begin{lmm}\label{sec:proof-theor-refs-2}
For every $\eta>0$, one has
$$\lim_{\eps\to
  0}\sup_{n\to\infty}\bP_n\left(\dgh(\bt,\bt[\lf\eps
  n\rf])>\eta\sqrt{n}\right)=0\, .$$
\end{lmm}

\proof Conditionally on its skeleton ${\rm Sk}_{\lf \eps
  n\rf}=(\ms,\bx,\by)$, and conditionally on the set $G_{\lf \eps
  n\rf}$ of probability going to $1$, the whole tree is reconstructed
by grafting trees all with sizes $\leq \lf \eps n\rf$ to the tree
$\Pi_0^{-1}(\ms,\bx)$, as in Lemma \ref{sec:case-trees-2}. Let
$r_1,r_2,\ldots$ be these sizes, in decreasing order, and
$\bt_1,\bt_2,\ldots$ be the associated trees (with some arbitrary
convention for ties). Of course, $\sum_i r_i=n$. Conditionally on the
sizes, the trees are random elements, respectively with distribution
$\bP_{r_i}$, by Lemma \ref{sec:case-trees-2}.

At this point, we rely on the following particular case of
\cite[Theorem 5]{BrFl08}, stating that, if $H:\bt\mapsto H(\bt)$
denotes the height (i.e.\ maximal height of elements of any $\t\in
\bt$) of the canonical random variable, then
$$\bE_n[H^4]\leq C n^2\, ,\qquad n\geq 1\, ,$$
where $\bE_n$ denotes expectation with respect to $\bP_n$, and $C$ is
some constant in $(0,\infty)$. Consequently, by the Markov inequality,
$$\bP_n(H\geq x\sqrt{n})\leq Cx^{-4}\, ,\qquad x>0\, .$$
Now, using this fact after conditioning on the sizes $r_i,i\geq 1$,
\begin{eqnarray*}
  \bP_n\left(\dgh(\bt,\bt[\lf\eps
    n\rf])>\eta\sqrt{n}\right) 
&\leq & \bP_n(\exists\, i\geq 0:H(\bt_i)\geq \eta\sqrt{n})\\
  &\leq & \bE_n\left[\sum_i P_{r_i}(H\geq \eta\sqrt{n})\right]\\
&\leq & C\eta^{-2}\bE_n\left[\sum_i (r_i/n)^2\right]\\
&\leq &C\eta^{-2}\eps\, .
\end{eqnarray*}
In the last step, we wrote $(r_i/n)^2\leq \eps(r_i/n)$ and summed over
$i$. Hence the result.  \cq

\medskip

To conclude the proof of Theorem \ref{sec:bc1}, simply apply the
following and last lemma, with $X_n=\bc(2n)^{-1/2}\bt$ under
$\bP_n$, $X_n^{(\eps)}=\bc(2n)^{-1/2}\bt[\lf \eps n\rf]$ under
$\bP_n$, $X^{(\eps)}=\theta(\zeta(2\ee,\eps))$ and $X=\TT_{2\ee}$. The
hypotheses are enforced, by Lemma \ref{sec:proof-theor-refs-2}, and
since $X_n^{(\eps)}\to X^{(\eps)}$ in distribution (Proposition
\ref{sec:interpr-terms-psims-1}) and $X^{(\eps)}\to X$ a.s.\ as a
consequence of Lemma \ref{sec:excurs-funct-trimm-2}.

\begin{lmm}\label{sec:proof-theor-refs-4}
  On some probability space $(\Omega,\FF,P)$, let
  $X_n,X,X_n^{(\eps)},X^{(\eps)},n\geq 1,\eps\in(0,1)$ be random
  variables with values in some metric space $(Z,d)$ such that  
$X_n^{(\eps)}\to X^{(\eps)}$ in distribution as $n\to\infty$,
  $X^{(\eps)}\to X$ in distribution as $\eps\to 0$, and
$$\lim_{\eps\to 0}\sup_{n\geq 1}P(d(X_n,X_n^{(\eps)})> \eta)=0\,.$$
Then $X_n\to X$ in distribution as $n\to\infty$.
\end{lmm}

\proof
Let $F$ be closed in $Z$, then 
$$P(X_n\in F)\leq P(X_n\in F,d(X_n,X_n^{\eps})\leq
\eta)+\delta_\eta^{(\eps)}\, ,$$ where $\delta_\eta^{(\eps)}=\sup_n
P(d(X_n,X_n^{(\eps)})> \eta)$. Thus $P(X_n\in F)\leq P(X_n^{(\eps)}\in
F_\eta)+\delta_\eta^{(\eps)}$, where $F_\eta=\{z\in Z:d(z,F)\leq
\eta\}$ is closed. Since $X_n^{(\eps)}\to X^{(\eps)}$ in distribution
as $n\to\infty$, it holds that 
$$  \limsup_n P(X_n\in F)\leq P(X^{(\eps)}\in F_\eta)+
  \delta_\eta^{(\eps)}\, ,$$
and taking the limsup as $\eps\to 0$, appealing to the convergence in
distribution $X^{(\eps)}\to X$, and finally, letting $\eta\to 0$, this
gives $\limsup_n P(X_n\in F)\leq P(X\in F)$
as wanted.  \cq

\section{Final comments}\label{sec:final-comments}

To our knowledge, the present work is the first to deal with the full
scaling limit picture of a random non-labeled, non-plane tree. As was
predictable from the known enumeration results, and conjectured by
Aldous, the CRT is the object that arises as the limit of uniform
random binary rooted trees. 

In this paper, we have purposely restricted our attention to binary
trees in order to make the combinatorial arguments as simple as
possible when considering trimming of trees. However, it is to be
expected that the kind of methods used to prove statements like
Proposition \ref{sec:interpr-terms-psims-1} will still work for more
general families of unordered trees, for instance trees which are `at
most $m$-ary' for some $2\leq m\leq \infty$, as considered in Otter's
work \cite{otter48} (the case $m=\infty$ is sometimes referred to as
{\em P\'olya trees}). What is to be expected is that after a trimming
by a mass $\lf \eps n\rf$ a uniform tree with $n$ vertices, the tree
that one obtains will be binary with a probability close to $1$,
because the CRT is itself binary.  In view of the probabilistic
approach taken in Section \ref{sec:case-trees-3}, it seems also
possible to adapt this method to more general families of trees, i.e.\
of trees satisfying an enumeration scheme similar to (\ref{eq:3}).

We have not gone through the details of this approach in full
generality, though, for several reasons. One of these is that it
needs a result similar to that of Broutin and Flajolet \cite{BrFl08}
for more general trees. This might appear as a minor issue, as
combinatorial methods are probably robust enough to be generalizable
to more general situations. For instance, while completing the present
work, we became aware of the paper \cite{drmgit09} by Drmota and
Gittenberger, in which moments estimates of the kind of \cite{BrFl08}
are obtained for the case of P\'olya trees. However, one should note
that the results of \cite{BrFl08} give much more than what is actually
needed to prove Lemma \ref{sec:proof-theor-refs-2}, i.e.\ precise
moment estimates, or the derivation of the exact scaling limit of the
height. From the probabilistic viewpoint, it would be very interesting
to be able to decide in a slicker way whether this `tightness'
property of Lemma \ref{sec:proof-theor-refs-2} holds in a general
setting. Note in particular that a derivation of the scaling limit for
the total height of the tree follows from Theorem \ref{sec:bc1}:
$$\bP_n(\bc \, H\geq x\sqrt{2n})\build\longrightarrow_{n\to\infty}^{}
P(2\max \ee\geq x)=2\sum_{k\geq 1}(k^2x^2-1)e^{-k^2x^2/2},~~~ \textrm{ for } x>0$$ (for the
last inequality, see e.g.\ \cite{DuIgDo}), which is \cite[Theorem
1]{BrFl08}.

\bibliographystyle{abbrv} \bibliography{biblio}

\end{document}